\documentclass[a4paper]{amsart} 

\usepackage{amsfonts}
\usepackage{amsmath}  
\usepackage{amscd}
\usepackage{graphicx}
\usepackage{psfrag}

\usepackage{amssymb}
\usepackage{latexsym}
\usepackage{framed, color}
\usepackage{mathabx}

\newtheorem{theorem}{\sc Theorem}[section]
\newtheorem{corollary}[theorem]{\sc Corollary}
\newtheorem{lemma}[theorem]{\sc Lemma}
\newtheorem{proposition}[theorem]{\sc Proposition}
\theoremstyle{definition}
\newtheorem{dfn}[theorem]{\sc Definition}

\theoremstyle{remark}
\newtheorem{exam}[theorem]{\sc Example}
\newtheorem{rmk}[theorem]{\sc Remark}
\newtheorem{prob}[theorem]{\sc Problem}

\makeatletter
 
 \@addtoreset{equation}{section}
\makeatother

\newcommand{\C}{\mathbf{C}}
\newcommand{\R}{\mathbf{R}}
\newcommand{\Q}{\mathbf{Q}}
\newcommand{\Z}{\mathbf{Z}}

\newcommand{\Int}{\mathop{\mathrm{Int}}\nolimits}
\newcommand{\Ker}{\mathop{\mathrm{Ker}}\nolimits}

\newcommand{\Image}{\mathop{\mathrm{Im}}\nolimits}

\renewcommand{\tilde}{\widetilde}
\renewcommand{\hat}{\widehat}
\renewcommand{\setminus}{\smallsetminus}
\def\spmapright#1{\smash{%
 \mathop{\hbox to 1.2cm{\rightarrowfill}}
  \limits^{#1}}}

\title{Cobordism of
algebraic knots defined by Brieskorn polynomials, II}
\author{Vincent Blanl\oe il and Osamu Saeki} 

\address{V.~Blanl\oe il: UFR Math\'{e}matique et Informatique, IRMA,
Universit\'{e} de Strasbourg,
7, rue Ren\'{e} Descartes,
F-67084 Strasbourg, France
}

\address{O.~Saeki: Institute of Mathematics for Industry,
Kyushu University, Motooka 744, Nishi-ku, Fukuoka 819-0395,
Japan}
\date{\today}

\begin{document}

\begin{abstract} 
In our previous paper, we obtained
several results concerning cobordisms of
algebraic knots associated with Brieskorn polynomials:
for example, under certain conditions, we showed that the
exponents are cobordism invariants.
In this paper, we further obtain new results
concerning the Fox--Milnor type relations, 
decomposition of the algebraic
cobordism class of an algebraic knot associated with
a Brieskorn polynomial that has a null-cobordant 
factor over the field of rational numbers,
and cyclic suspensions of knots.
We also show that a certain infinite
family of spherical algebraic
knots associated with Brieskorn polynomials are linearly
independent in the knot cobordism group.
\end{abstract}

\maketitle 

\section{Introduction}\label{section1}

\begin{figure}[tb]\label{fig-sing}
\centering
 \begin{picture}(120,120)(0,0)
  \put(50,50){\circle*{3}}
  \put(80,72){\circle*{3}}
  \put(80,28){\circle*{3}}
  \put(95,110){\makebox{$V_f$}}
  \put(0,90){\makebox{$S^{2n+1}_{\varepsilon}$}}
  \put(110,48){\makebox{$K_f$}}
  \put(47,55){\makebox{$\mathbf{0}$}}
  \put(105,60){\vector(-2,1){23}}
  \put(105,40){\vector(-2,-1){23}}
  \thicklines
  \qbezier(50,50)(80,65)(110,110)
  \qbezier(50,50 )(80,35)(110,-10)
  \thinlines
  \qbezier(80,72)(72,90)(50,90)
  \qbezier(50,90)(10,90)(10,50)
  \qbezier(10,50)(10,10)(50,10)
  \qbezier(50,10)(72,10)(80,28)
  \qbezier(80,28)(87,50)(80,72)
  \put(80,72){\circle*{3}}
  \put(80,28){\circle*{3}} 
 \end{picture}
\caption{The algebraic knot $K_f$ associated with the 
singularity at $\mathbf{0}$ of a germ $f$}

\end{figure}
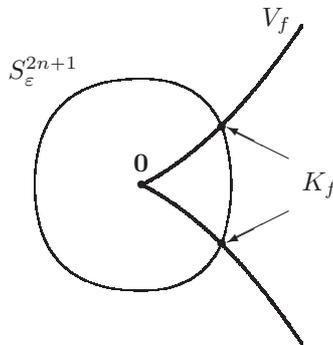

Let $f : (\C^{n+1}, \mathbf{0}) \to (\C, 0)$, $n \geq 1$, 
be a holomorphic
function germ with an isolated critical point at the origin.
For a sufficiently small positive real number $\varepsilon > 0$,
set $K_f = S^{2n+1}_\varepsilon \cap V_f$,
where $V_f = f^{-1}(0)$ is the complex hypersurface
in $\C^{n+1}$ with an isolated singularity at the origin and
$S^{2n+1}_\varepsilon$ is the sphere of radius $\varepsilon$
centered at the origin in $\C^{n+1}$ (see Fig.~\ref{fig-sing}).
It is known that $K_f$ is an $(n-2)$--connected, oriented
$(2n-1)$--dimensional submanifold of $S^{2n+1}_\varepsilon
= S^{2n+1}$, that its complement fibers over the circle $S^1$,
and that the isotopy class of $K_f$ in $S^{2n+1}$
is independent of the choice of $\varepsilon$ as long as
it is sufficiently small
(see \cite{Milnor}). Note also that 
the embedded topology of $V_f \subset \C^{n+1}$
around the origin determines and is determined by
the (oriented) isotopy class of $K_f \subset S^{2n+1}$
(see \cite{Sa89}).
We call $K_f$ the \emph{algebraic knot}
associated with $f$.
In this paper, a \emph{knot} (or a \emph{$(2n-1)$--knot})
refers to (the isotopy
class of) an $(n-2)$--connected, oriented $(2n-1)$--dimensional
submanifold of $S^{2n+1}$. (Here, when $n=1$, a submanifold
is $(-1)$--connected if it is nonempty.)
A $(2n-1)$--knot $K$ is \emph{spherical} if it is homeomorphic
to the sphere $S^{2n-1}$.

In this paper, we consider Brieskorn polynomials
\begin{equation}
f(z_1, z_2, \ldots, z_{n+1}) =
z_1^{a_1} + z_2^{a_2} + \cdots + z_{n+1}^{a_{n+1}}
\label{eq:Br}
\end{equation}
with exponents $a_i \geq 2$, $1 \leq i \leq n+1$, and their
associated algebraic knots $K_f$ \cite{Brieskorn}.
We especially focus on the study of their properties
concerning cobordisms.
Two knots $K_0$ and $K_1$ in $S^{2n+1}$
are said to be \emph{cobordant}
if there exists a properly embedded oriented submanifold $X$,
abstractly diffeomorphic to $K_0 \times [0, 1]$,
of $S^{2n+1} \times [0, 1]$ such that
$X \cap (S^{2n+1} \times \{0\}) = K_0$, and
$X \cap (S^{2n+1} \times \{1\}) = -K_1^!$,
where $-K_1^!$ is the mirror image of $K_1$ with 
reversed orientation (see \cite{BM, BS}).
We say that a $(2n-1)$--knot $K$ in $S^{2n+1}$
is \emph{null-cobordant} if it bounds a smoothly
embedded $2n$--dimensional disk in $D^{2n+2}$
(note that $\partial D^{2n+2} = S^{2n+1}$).

In our previous paper \cite{BS2}, we obtained
several results concerning cobordisms of
algebraic knots associated with Brieskorn polynomials:
for example, under certain conditions, we showed that the
exponents are cobordism invariants.
In this paper, we further obtain new results
concerning the Fox--Milnor type relations for Alexander polynomials, 
decomposition of the algebraic
cobordism class of a spherical algebraic knot associated with
a Brieskorn polynomial that has a null-cobordant 
factor over the field of rational numbers,
and cyclic suspensions of knots.

The present paper is organized as follows.
In \S\ref{section2}, we recall several basic definitions
and properties concerning invariants and cobordisms of 
algebraic knots
such as Alexander polynomials and Seifert forms.

In \S\ref{section3}, we focus on the Fox--Milnor type relations
for Alexander polynomials \cite{F-M0, F-M} and
give a complete characterization of Brieskorn
polynomials that give algebraic knots whose Alexander polynomials
satisfy the Fox--Milnor type relation in terms of their exponents.
As a consequence, we show that an algebraic knot
associated with a Brieskorn polynomial is never null-cobordant:
moreover, it turns out that a spherical algebraic knot associated
with a Brieskorn polynomial always has infinite order in the
knot cobordism group. 
In fact, Michel \cite{Michel} has shown that
such a result holds for algebraic knots in general:
our proof shows that an argument based on 
the Fox--Milnor type relation
serves well for deducing such results
at least for Brieskorn polynomials.

In \S\ref{section3.5}, we consider the linear independence
of a family of spherical algebraic knots associated with 
certain Brieskorn
polynomials in the knot cobordism group.
In fact, Litherland \cite{Litherland} has shown that the spherical
algebraic knots in $S^3$ associated with Brieskorn 
polynomials of two variables (in fact, such knots are
torus knots) are linearly independent
in the $1$--dimensional knot cobordism group by using 
a certain signature invariant. We will use the same idea
to prove a similar linear independence result for higher dimensions.

In \S\ref{section4}, we consider the group of algebraic cobordism
classes of spherical knots which has been introduced and
studied by Levine \cite{L1, L2} (see also \cite{Collins}). 
We give an explicit
example of a spherical algebraic knot
associated with a Brieskorn polynomial
such that its algebraic cobordism class has a
decomposition into those corresponding to
the irreducible factors of its Alexander polynomial
over the field of rational numbers and that
one of them is algebraically null-cobordant.
This shows that cobordant spherical algebraic knots
associated with Brieskorn polynomials 
may not share the same irreducible factors of their
Alexander polynomials, and therefore
the study of cobordism classes of
algebraic knots associated with Brieskorn polynomials
might be more complicated than is expected.

Finally in \S\ref{section5}, we consider cyclic suspensions of
knots \cite{Kauffman, N} and study its relationship to
the cobordisms. Note that the algebraic knot
associated with a polynomial of the form
$f(z_1, z_2, \ldots, z_{n+1}) + z_{n+2}^d$
is the $d$--fold cyclic suspension of the algebraic
knot associated with $f$.
We will see that the cyclic suspension
of knots often behaves very badly with respect to cobordisms.
For example we show that certain cyclic suspensions
of the algebraic knots
constructed by Du Bois--Michel in \cite{dBM}, which
are cobordant to each other, are not diffeomorphic
and are not cobordant.

Throughout the paper, all manifolds and maps
between them are smooth of class $C^\infty$.
The symbol ``$\cong$'' means a diffeomorphism
between manifolds or an appropriate isomorphism between
algebraic objects.

\section{Preliminaries}\label{section2}

This paper is a sequel of our previous paper \cite{BS2}.
The reader is expected to be familiar with the
notions and basic results explained in that paper,
although we will repeat them when necessary
in the present paper.

Let $K$ be a $(2n-1)$--knot in $S^{2n+1}$.
Suppose that there exists a locally trivial
fibration $\varphi : S^{2n+1} \setminus K \to S^1$.
We also assume that there is a trivialization 
$\tau : N(K) \to K \times D^2$ of the normal disk
bundle neighborhood $N(K)$ of $K$ in $S^{2n+1}$ such that
the composition
$$N(K) \setminus K \spmapright{\tau|_{N(K) \setminus K}}
K \times (D^2 \setminus \{0\}) \spmapright{pr_2} 
D^2 \setminus \{0\} \spmapright{r} S^1$$
coincides with $\varphi|_{N(K) \setminus K}$, where
$pr_2$ is the projection to the second factor and
$r$ is the radial projection.
Then, we say that $K$ is a \emph{fibered knot}.
We call the closure $F$ of a fiber of $\varphi$ a
\emph{fiber}. Note that it is a $2n$--dimensional
compact oriented submanifold of $S^{2n+1}$ whose boundary
coincides with $K$.
A $(2n-1)$--dimensional fibered knot $K$ is
\emph{simple} if it is $(n-2)$--connected and
$F$ is $(n-1)$--connected.
(Here, for $n=1$, a manifold is $(-1)$--connected if it
is nonempty.)
In this case, $F$ is homotopy equivalent to a bouquet
of $n$--dimensional spheres (for example, see \cite[Theorem~6.5]{Milnor}).
Note that an algebraic knot associated with a holomorphic
function germ $f : (\C^{n+1}, \mathbf{0})
\to (\C, 0)$ with an isolated critical point at the origin
is a simple fibered knot \cite{Milnor}. In this case,
a fiber of such an algebraic knot is called a \emph{Milnor
fiber} for $f$.

Let $\psi : F \to F$ be a \emph{geometric monodromy}
of the fibration $\varphi$; i.e., it is a diffeomorphism
which is constructed by integrating an appropriate
horizontal vector field on $S^{2n+1} \setminus K$ with respect
to $\varphi$ and
which is the identity on the boundary. 
In other words, $S^{2n+1} \setminus \Int{N(K)}$
is diffeomorphic to the manifold
$$F \times [0, 1]/(x, 1) \sim (\psi(x), 0), \quad x \in F,$$
obtained by identifying $F \times \{1\}$ and $F \times \{0\}$ by $\psi$.
It is known that the geometric monodromy is
well-defined up to isotopy.
Either of the isomorphisms
$$\psi_* : H_n(F; \Z) \to H_n(F; \Z) \quad \mbox{or}
\quad \psi^* : H^n(F; \Z) \to H^n(F; \Z)$$
is called the \emph{algebraic monodromy}.
Its characteristic polynomial $\Delta_K(t) \in \Z[t]$, 
which is well-defined for both of $\psi_*$ and $\psi^*$,
is often called the \emph{Alexander polynomial} of $K$.
When $K$ is an algebraic knot associated with
a holomorphic function germ $f$, we often denote $\Delta_{K_f}(t)$
by $\Delta_f(t)$.

Let us consider the multiplicative group $\C^*$
and its group ring $\Z \C^*$ over the integers.
For a monic polynomial $\Delta(t)$ with nonzero
constant term, we denote by
$\mathrm{divisor}\,\Delta$ the element
$$\sum m_\xi \langle \xi \rangle \in \Z \C^*,$$
where $\xi \in \C^*$ runs over all roots of 
$\Delta(t)$ and $m_\xi \in \Z$
is its multiplicity. 

We also set
$$\Lambda_a = \mathrm{divisor}\,(t^a -1)$$
for a positive integer $a$.
Note that the family of such elements $\{\Lambda_a\}_a$,
where $a$ runs over all positive integers, is linearly
independent over $\Z$.

Now, let us consider a Brieskorn polynomial as in (\ref{eq:Br}).
Then, by Brieskorn \cite{Brieskorn}, it is known that
\begin{equation}
\mathrm{divisor}\, \Delta_f =
\left(\Lambda_{a_1} - 1\right)
\left(\Lambda_{a_2} - 1\right) \cdots
\left(\Lambda_{a_{n+1}} - 1\right),
\label{eq:BR}
\end{equation}
where $1 \in \Z \C^*$ means $\Lambda_1 = 1\langle 1 \rangle$.
This implies that the roots of $\Delta_f(t)$
are all roots of unity and that $\Delta_f(t)$
is a product of cyclotomic polynomials.
In particular, each irreducible factor $\gamma(t)
\in \Z[t]$ of $\Delta_f(t)$ is \emph{symmetric}, i.e.,
we have
$\gamma(t) = \pm t^{\mathrm{deg}\,\gamma}\gamma(t^{-1})$,
where $\mathrm{deg}\,\gamma \in \Z$ is the degree
of the polynomial $\gamma$.

Let $K$ be a $(2n-1)$--knot.
We say that $K$ is \emph{spherical} if $K$ is
homeomorphic to the $(2n-1)$--dimensional sphere.
This means that $K$ may be an exotic sphere, which
is homeomorphic but not diffeomorphic to the standard sphere.
When $K$ is a simple fibered $(2n-1)$--knot with $n \neq 2$,
it is known that $K$ is spherical if and only if $\Delta_K(1)
= \pm 1$ (for example, see \cite[Theorem~8.5]{Milnor}).
For algebraic knots associated with a Brieskorn
polynomial, there is a characterization of spherical knots
due to Brieskorn \cite{Brieskorn} in terms of the exponents
(for details, see Theorem~\ref{thm:Br} and 
Remark~\ref{rem:Mumford} of the present paper).

\begin{figure}[tb]
\centering
\begin{picture}(100,100)(0,0)
\thicklines
\qbezier(10,50)(-30,90)(50,90)
\put(53,90){\vector(-1,0){0}}
\qbezier(50,90)(120,90)(108,71)
\qbezier(105,67)(98,57)(92,52)
\qbezier(50,50)(70,35)(88,48)
\qbezier(12,52)(30,65)(50,50)
\qbezier(8,48)(-30,10)(50,10)
\qbezier(50,10)(130,10)(90,50)
\qbezier(90,50)(70,65)(52,52)
\qbezier(48,48)(30,35)(10,50)
\put(15,100){\shortstack{$+$}}
\put(15,82){\vector(1,3){5}}
\thinlines
\qbezier(51,52)(50,65)(59,72)
\qbezier(61,73)(70,75)(78,72)
\qbezier(80,71)(90,67)(90,52)
\put(70,48){\oval(39,40)[b]}
\qbezier(30,43)(25,75)(50,75)
\qbezier(50,75)(72,75)(65,44)
\qbezier[25](30,43)(50,2)(64,41)
\qbezier(11,52)(12,77)(37,73)
\qbezier(40,72)(50,67)(49,52)
\put(30,48){\oval(39,40)[b]}
\qbezier(70,43)(70,75)(90,75)
\qbezier(90,75)(120,75)(100,41)
\qbezier[17](98,38)(70,15)(70,43)
\put(32,78){$a$}
\put(33,74){\vector(1,0){0}}
\put(50,78){$a_{+}$}
\put(53,75){\vector(1,0){0}}
\put(65,78){$b$}
\put(73,74){\vector(1,0){0}}
\put(85,78){$b_{+}$}
\put(93,75){\vector(1,0){0}}
\put(100,10){$F$}
\end{picture}
\caption{Computing a Seifert matrix for the trefoil knot}
\label{fig1}
\end{figure}
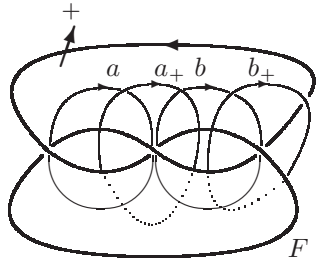

Let $K$ be a simple fibered $(2n-1)$--knot with fiber $F$.
We define the bilinear form
$\theta_K : H_n(F; \Z) \times H_n(F; \Z) \to \Z$ by
$\theta_K(\alpha, \beta) = \mathrm{lk}(a_+, b)$, where
$a$ and $b$ are $n$--cycles representing $\alpha$
and $\beta$, respectively, $a_+$ is the
$n$--cycle in $S^{2n+1}$ obtained by pushing $a$ 
into the positive normal direction of $F$, and $\mathrm{lk}$ 
denotes the
linking number of $n$--cycles in $S^{2n+1}$
(see Fig.~\ref{fig1}).
The bilinear form $\theta_K$ is called the \emph{Seifert form}
of $K$ and its representative matrix is called a \emph{Seifert matrix}.
It is known that a Seifert form is unimodular, i.e., the
determinant of the Seifert matrix $L_K$
is equal to $\pm 1$,
due to Alexander duality (see \cite{Durfee}).
Furthermore, it is also known that the Alexander
polynomial $\Delta_K(t)$ coincides with
$\pm \det (t L_K + (-1)^n L_K^T)$, where
$L_K^T$ denotes the transpose of $L_K$
(for example, see \cite{Sa00}).

It is known that for $n \geq 3$, there is a one-to-one
correspondence, through Seifert forms, 
between the set of isomorphism classes
of unimodular bilinear forms over the integers and the set
of isotopy classes of simple fibered $(2n-1)$--knots
\cite{Durfee, Kato}.

Recall that the set of cobordism classes of
spherical $(2n-1)$--knots forms an additive group
under the connected sum operation. This
is called the \emph{$(2n-1)$--dimensional knot
cobordism group} and is denoted by $C_{2n-1}$
(for example, see \cite{Collins, L1, L2}.)
The class of the trivial knot is the neutral
element, and the inverse of the cobordism class of a knot $K$ is
the class of $-K^{!}$.
Note that for $n > 1$, it is known that $C_{2n-1}
\cong \Z^\infty \oplus \Z_2^\infty \oplus \Z_4^\infty$
(see \cite{L1, L2}).

Let us now recall the algebraic cobordism group $G_\varepsilon$,
where $\varepsilon = (-1)^n$ (for details, see \cite{L1}).
We consider square integer matrices $A$ such
that $A + \varepsilon A^T$ is unimodular:
such a matrix is called an $\varepsilon$--matrix.
An $\varepsilon$--matrix $N$ is \emph{null-cobordant} if
$N$ is congruent to a matrix of the form
$$\begin{pmatrix} 0 & N_1 \\ N_2 & N_3
\end{pmatrix},$$
where $N_1, N_2$ and $N_3$ are square matrices of
the same size.
Two $\varepsilon$--matrices $A_1$ and $A_2$ are
\emph{cobordant} if $A_1 \oplus (-A_2)$ is null-cobordant.
Then, it is known that this defines an equivalence relation
for $\varepsilon$--matrices, and the set $G_\varepsilon$
of cobordism classes of $\varepsilon$--matrices forms
an abelian group under block sum ``$\oplus$''.
Then, by Levine \cite{L1}, it has been proved that
for $n > 2$, the knot cobordism group $C_{2n-1}$
is isomorphic to $G_\varepsilon$ with $\varepsilon = (-1)^n$.

Now, let us consider square matrices $B$ with entries
in the field $\Q$ of rational numbers. We say that $B$ is
\emph{admissible} if
$$(B - B^T)(B + B^T)$$
is nonsingular.
Then, the cobordism relation is also defined for admissible
matrices, and the set $G^\Q$ of cobordism classes of
admissible matrices again forms an abelian group
under block sum.
Furthermore, it is known that the natural inclusion
$G_\varepsilon \to G^\Q$ is a monomorphism.
Furthermore, a complete set of invariants for $G^\Q$
has been given by Levine \cite{L2}.

It is also known that  $G^{\Q}$ is isomorphic to
the group $G_{\Q}$
of cobordism classes of \emph{isometric structures} over $\Q$
(see \cite{Collins}).

Let us now consider the case of algebraic knots.
Let $f : (\C^{n+1}, \mathbf{0})
\to (\C, 0)$ be a holomorphic function germ
with an isolated critical point at the origin.
Let us define several notions concerning
the Seifert form $\theta_{K_f}$ of the algebraic knot
$K_f$ associated with $f$.

\begin{dfn}
Two bilinear forms $\theta_i : H_i \times H_i \to \Z$,
$i = 0, 1$, defined on free abelian groups $H_i$ of
finite ranks are said to be \emph{Witt equivalent}
(or \emph{cobordant})
if there exists a direct summand $M$ of $H_0 \oplus
H_1$ such that $(\theta_0 \oplus (-\theta_1))(x, y) = 0$
for all $x, y \in M$ and twice the rank of $M$
is equal to the rank of $H_0 \oplus H_1$.
In this case, $M$ is called a \emph{metabolizer}.

Furthermore, we say that $\theta_0$ and $\theta_1$
are \emph{Witt equivalent over the real numbers}
if there exists a vector subspace $M_{\R}$
of $(H_0 \otimes \R) \oplus (H_1 \otimes \R)$
such that $(\theta_0^\R \oplus (-\theta_1^\R))(x, y)
= 0$ for all $x, y \in M_{\R}$ and
$2\dim_\R M_\R = \dim_\R (H_0 \otimes \R)
+ \dim_\R (H_1 \otimes \R)$, where
$\theta_i^\R : (H_i \otimes \R) \times (H_i \otimes \R)
\to \R$ is the real bilinear form associated with $\theta_i$,
$i = 0, 1$.
\end{dfn}

Now, if $n$ is odd, let us consider 
$$\tilde{f}(z_1, z_2, \ldots, z_{n+1}, z_{n+2})
= f(z_1, z_2, \ldots, z_{n+1}) + z_{n+2}^2.$$
It is known that the Seifert form $\theta_{K_f}$ for $K_f$
coincides with $\theta_{f_{\tilde{f}}}$ for $K_{\tilde{f}}$
(for example, see \cite{Saka}). In the following, we may
assume that $n$ is even.
We have the decomposition
$$
H^n(F_f; \C) = \oplus_\lambda H^n(F_f; \C)_\lambda,
$$
where $F_f$ is the Milnor fiber for $f$,
$\lambda$ runs over all
the complex roots of the Alexander polynomial $\Delta_f(t)$, and
$H^n(F_f; \C)_\lambda$ is the eigenspace of the algebraic
monodromy $H^n(F_f; \C) \to H^n(F_f; \C)$ with respect
to the complex coefficients
corresponding to the eigenvalue $\lambda$.
It is known that the intersection form $S_f$ of $F_f$ 
on $H^n(F_f; \C)$ is given by
$S_f = L_f + L_f^T$ (refer to (\ref{eq:matrix}) of
\S\ref{section4}). Furthermore, this
decomposes as the orthogonal direct sum of 
$(S_f)|_{H^n(F_f; \C)_\lambda}$.
Let 
$\mu(f)_\lambda^+$ (resp.\ $\mu(f)_\lambda^-$) denote the 
number of positive (resp.\ negative) eigenvalues 
of $(S_f)|_{H^n(F_f; \C)_\lambda}$. 
Then, the integer
$$
\sigma_\lambda(f) = \mu(f)_\lambda^+ - \mu(f)_\lambda^-
$$
is called the \emph{equivariant signature} of 
$K_f$ with respect to $\lambda$ (for details,
see \cite{Neumann, SSS}).
Note that by the same construction, the equivariant
signatures are defined for fibered $(2n-1)$--knots
in general. For spherical $(2n-1)$--knots, by
appropriately defining an \emph{isometric structure} by
using a Seifert form, one can also define their
equivariant signatures (see \cite{L2}).

Then, it is known that the equivariant signatures
are integer-valued cobordism invariants. Furthermore, 
they are additive under connected sum.
In fact, it has been known that 
if two knots are cobordant, then their Seifert forms
are Witt equivalent over the real numbers, and that
two Seifert forms are Witt
equivalent over the real numbers if and only if
all the equivariant signatures coincide.
For details, the reader is referred to \cite{BS2}.

\section{Fox--Milnor type relation}\label{section3}

Let $\Delta_f(t)$ and $\Delta_g(t)$ denote
the Alexander polynomials for the algebraic knots
$K_f$ and $K_g$ associated with $f$ and $g$,
respectively, where $f$ and $g : 
(\C^{n+1}, \mathbf{0}) \to (\C, 0)$, $n \geq 1$, 
are holomorphic
function germs with an isolated critical point at the origin.
We say that the Alexander polynomials
satisfy the \emph{Fox--Milnor type relation} if
there exists a polynomial $\gamma(t)$ with integer coefficients
such 
that $\Delta_f(t)\Delta_g(t) = \pm t^{\mathrm{deg}\, \gamma}
\gamma(t) \gamma(t^{-1})$ (\cite{F-M0, F-M}).
It is known that if
$K_f$ and $K_g$ are cobordant, then their Alexander
polynomials satisfy the Fox--Milnor type relation
(for details, see \cite{F-M0, F-M, BM}.
See also \cite{BS}, for example). 

We have the following characterization
of Alexander polynomials which satisfy
the Fox--Milnor type relation.
In the following, for two elements
$d_1$ and $d_2 \in \Z\C^*$, we write
$d_1 \equiv d_2 \pmod{2}$ if there
exists an element $d_3 \in \Z\C^*$
such that $d_1 - d_2 = 2 d_3$
holds in $\Z\C^*$.

\begin{lemma}\label{rem:square}
For algebraic knots $K_f$ and $K_g$ as above,
the following three are equivalent to each other.

$(1)$ The Alexander polynomials
$\Delta_f(t)$ and $\Delta_g(t)$ satisfy the
Fox--Milnor type relation.

$(2)$ We have $\Delta_f(t) \Delta_g(t) = \gamma(t)^2$ for 
some $\gamma(t) \in \Z[t]$.

$(3)$ We have $\mathrm{divisor}\, \Delta_f(t) \equiv
\mathrm{divisor}\, \Delta_g(t) \pmod{2}$.
\end{lemma}

\begin{proof}
By \cite{Brieskorn70}, the
Alexander polynomials $\Delta_f(t)$ and $\Delta_g(t)$
are products of cyclotomic polynomials; in particular,
each of their irreducible factors is symmetric.
Therefore, their Alexander polynomials satisfy the
Fox--Milnor type relation if and only if $\Delta_f(t)\Delta_g(t)
= \gamma(t)^2$ for some $\gamma(t) \in \Z[t]$.

On the other hand, if $\Delta_f(t)\Delta_g(t)
= \gamma(t)^2$ for some $\gamma(t) \in \Z[t]$,
then obviously the congruence in $(3)$ holds.
Conversely, if the congruence in $(3)$ holds, then
as $\Delta_f(t)\Delta_g(t)$ is a product of
cyclotomic polynomials, each irreducible factor
appears an even number of times, so that $(2)$
holds. This completes the proof.
\end{proof}

Let
$$f(z) = z_1^{a_1} + z_2^{a_2} + \cdots + z_{n+1}^{a_{n+1}}$$
be a Brieskorn polynomial with $a_j \geq 2$ for all $j$. 
Set $E_f = \{a_1, a_2, \ldots, a_{n+1}\}$, which
may contain the same integer multiple times and is considered
to be a \emph{multi-set}.

\begin{dfn}\label{dfn:Ebar}
From $E_f$, we construct
the (non multi-)subset $\overline{E}_f \subset E_f$ 
by the successive procedure as follows.
\begin{enumerate}
\item Take off all those even integers which appear an even number of times. \label{(1)}
\item Take off the multiple elements except for one in such
a way that we get a non multi-set.
\label{(2)}
\item Take off $a_j$ if it is an integer multiple
of an odd $a_k$ with $k \neq j$.
\label{(3)}
\end{enumerate}
We call the set $\overline{E}_f$ thus obtained the \emph{essential
exponent set} of $f$. Note that $\overline{E}_f$ can be
empty.
\end{dfn}

The first main result of this paper is the following.

\begin{theorem}\label{thm:FM}
Let
$$f(z) = z_1^{a_1} + z_2^{a_2} + \cdots + z_{n+1}^{a_{n+1}}
\mbox{ and } 
g(z) = z_1^{b_1} + z_2^{b_2} + \cdots + z_{n+1}^{b_{n+1}}
$$
be Brieskorn polynomials with $a_j \geq 2$ and $b_j \geq 2$
for all $j$. 
Then, the Alexander polynomials $\Delta_f(t)$
and $\Delta_g(t)$ satisfy the Fox--Milnor type
relation if and only if their essential exponent sets
coincide, i.e.\ $\overline{E}_f = \overline{E}_g$.
\end{theorem}

\begin{exam}
For example, consider
$$f(z) = z_1^{3} + z_2^{4} + z_3^{4} + z_4^{6} + z_5^{9}
\mbox{ and }
g(z) = z_1^2 + z_2^2 + z_3^3 + z_4^3+z_5^{12}.
$$
Then, we have
$$E_f = \{3, 4, 4, 6, 9\} \mbox{ and }
E_g = \{2, 2, 3, 3, 12\}.$$
In the process of Definition~\ref{dfn:Ebar},
after (\ref{(1)}), we get the multi-sets 
$\{3, 6, 9\}$ and $\{3, 3, 12\}$
for $f$ and $g$, respectively.
After (\ref{(2)}), we get the sets $\{3, 6, 9\}$ and $\{3, 12\}$.
Finally, after (\ref{(3)}), we get the sets $\{3\}$ and $\{3\}$.
Hence, we get $\overline{E}_f = \overline{E}_g = \{3\}$
and $\Delta_f(t)$ and $\Delta_g(t)$ satisfy the Fox--Milnor
type relation. In fact, by the formula (\ref{eq:BR})
with the help of Lemma~\ref{lemma:basic} below, we have
\begin{eqnarray*}
\mathrm{divisor}\,\Delta_f(t) & = & 4 \Lambda_{12}
- \Lambda_3 - 1, \\
\mathrm{divisor}\,\Delta_g(t) & = &
24 \Lambda_{36} + 6\Lambda_{18}
-6\Lambda_{12} -2\Lambda_9 -2\Lambda_6
-2\Lambda_4+\Lambda_3-1,
\end{eqnarray*}
so we can verify that
$\Delta_f(t)$ and $\Delta_g(t)$ satisfy the Fox--Milnor type
relation by virtue of Lemma~\ref{rem:square}.

Note that by the signature formula due to Brieskorn \cite{Brieskorn},
we see that the signatures of the $8$--dimensional
Milnor fibers for $f$ and $g$ are
equal to $274$ and $30$, respectively. Thus, $K_f$ and $K_g$
are not cobordant, since the signature of a fiber of a fibered knot
is a cobordism invariant, which can be proved
by using the fact that cobordant fibered knots have
algebraically cobordant Seifert forms \cite{BM}.
Nevertheless, their Alexander polynomials satisfy
the Fox--Milnor type relation.

On the other hand, for
$$h(z) =  z_1^{3} + z_2^{4} + z_3^{4} + z_4^{6} + z_5^{8},$$
we have $\overline{E}_h = \{3, 8\}$, so $\Delta_f(t)$ (or $\Delta_g(t)$)
and $\Delta_h(t)$ do not satisfy the Fox--Milnor type
relation. In fact, we have
$$\mathrm{divisor}\,\Delta_h(t) = 27\Lambda_{24}
-6\Lambda_{12} + 9\Lambda_8 -2\Lambda_6
-2\Lambda_4+\Lambda_3-1,$$
which verifies the above assertion by virtue of
Lemma~\ref{rem:square}.
\end{exam}

\bigskip

In order to prove Theorem~\ref{thm:FM}, let us prepare
some preliminary lemmas.
Recall that we have
$$\mathrm{divisor}\,\Delta_f(t) = \prod_{i=1}^{n+1}
(\Lambda_{a_i} - 1)\quad  \mbox{ and } \quad
\mathrm{divisor}\,\Delta_g(t) = \prod_{i=1}^{n+1}
(\Lambda_{b_i} - 1) 
$$
by \cite{Brieskorn}
and that $\Delta_f(t)$ and $\Delta_g(t)$ satisfy the
Fox--Milnor type relation if and only if
$$\mathrm{divisor}\,\Delta_f(t) \equiv
\mathrm{divisor}\,\Delta_g(t)  \pmod{2}$$
by Lemma~\ref{rem:square}.

In the following, we will use
the basic formula as follows
(for example, see \cite{MO}).
For completeness, we give a proof below.

\begin{lemma}\label{lemma:basic}
For positive integers $a$ and $b$, we have
\begin{equation}
\Lambda_a \Lambda_b = (a, b)\Lambda_{[a, b]},
\label{eq:basic}
\end{equation}
where $(a, b)$
denotes the greatest common divisor of $a$ and $b$, and
$[a, b]$ denotes the least common multiple of $a$ and $b$.
\end{lemma}

\begin{proof}
For $a=1$ or $b=1$, the conclusion is obvious. So, we may
assume $a>1$ and $b>1$.
Set $d = (a, b)$, $a' = a/d$ and $b' = b/d$, which
are positive integers. Note that then
$(a', b') = 1$ and $[a, b] = da'b'$. 
Let $\xi$ be a primitive $(da'b')$--th
root of unity. Then, we have
\begin{equation}
\Lambda_a \Lambda_b 
= \left(\sum_{k=0}^{a-1}\langle \xi^{b'k} \rangle \right)
\left(\sum_{\ell=0}^{b-1}\langle \xi^{a'\ell} \rangle \right)
= \sum_{k=0}^{a-1} \sum_{\ell=0}^{b-1}
\langle \xi^{b'k + a'\ell} \rangle.
\label{eq:sum}
\end{equation}
Now, we see easily that
$\langle \xi^{b'k + a'\ell} \rangle = \langle 1 \rangle$
for 
$$(k, \ell) = (0, 0), (a', b-b'), (2a', b-2b'),
\ldots, ((d-1)a', b-(d-1)b').$$
Therefore $\langle 1 \rangle$ appears at least
$d$ times in the rightmost summation in (\ref{eq:sum}). 
Similarly, given a pair of integers $(\alpha, \beta)$
with $0 \leq \alpha \leq a-1$ and $0 \leq \beta \leq b-1$,
we have 
$\langle \xi^{b'k + a'\ell} \rangle = \langle 
\xi^{b'\alpha + a'\beta} \rangle$
for 
\begin{eqnarray}
(k, \ell)  & \equiv & (\alpha, \beta), (\alpha+a', \beta-b'),
(\alpha+2a', \beta-2b'), 
\label{eq:list}
\\
& & \quad \ldots, (\alpha + (d-1)a',
\beta - (d-1)b') \pmod{(a, b)}
\nonumber
\end{eqnarray}
with $0 \leq k \leq a-1$ and $0 \leq \ell \leq b-1$,
where for integers $r, r', s$ and $s'$, we write
$(r, s) \equiv (r', s') \pmod{(a, b)}$
when $r \equiv r' \pmod{a}$ and
$s \equiv s' \pmod{b}$. Note that
the integer pairs in (\ref{eq:list}) are all
distinct modulo $(a, b)$. Therefore, $\langle 
\xi^{b'\alpha + a'\beta} \rangle$ also appears 
at least
$d$ times in the rightmost summation in (\ref{eq:sum}). 

When $(\alpha, \beta)$ runs over all
integer pairs
with $0 \leq \alpha \leq a-1$ and $0 \leq \beta \leq b-1$,
we see that $\langle 
\xi^{b'\alpha + a'\beta} \rangle$ runs over all
$\langle \xi^{\gamma} \rangle$ with
$0 \leq \gamma \leq da'b'-1$, since $a'$ and $b'$
are relatively prime. Hence, we see that each
$\langle \xi^{\gamma} \rangle$ appears exactly $d$ times.
Therefore, we have the desired conclusion.
This completes the proof.
\end{proof}

Using the basic formula in Lemma~\ref{lemma:basic},
we can show the following.

\begin{lemma}\label{lemma:FM}
For positive integers $a, b$ and $m$, we have the following.
\begin{enumerate}
\item If $a$ is even, then we have
$$(\Lambda_a-1)^m \equiv
\begin{cases} 1 \pmod{2}, & \mbox{$m$: even}, \\
\Lambda_a-1 \pmod{2}, & \mbox{$m$: odd}.
\end{cases}
$$
\item If $a$ is odd, then we have
$$(\Lambda_a-1)^m \equiv \Lambda_a -1 \pmod{2}$$
for all $m$.
\item If $a$ is odd, then we have
$$(\Lambda_a-1)(\Lambda_{ab} -1) \equiv
\Lambda_a-1 \pmod{2}.
$$
\item If $a_j$, $j = 1, 2, \ldots, m$, are positive even integers, 
then we have
$$\prod_{j=1}^m (\Lambda_{a_j} - 1) \equiv 
\sum_{j=1}^m \Lambda_{a_j} - 1 \pmod{2}.
$$
\item If $a$ is even and $b$ is odd, then we have
$$\Lambda_a (\Lambda_b-1) \equiv \Lambda_{[a, b]} - \Lambda_a
\pmod{2}.$$ 
\end{enumerate}
\end{lemma}

\begin{proof}
(1) Let us prove the assertion by induction on $m$.
For $m=1$, it is obvious. For $m=2$, by
Lemma~\ref{lemma:basic} we have
$$(\Lambda_a-1)^2 = a\Lambda_a - 2 \Lambda_a + 1
\equiv 1 \pmod{2},$$
since $a$ is even. So, the result holds.
Then, suppose $m \geq 3$ and that the assertion holds
for $m-1$. If $m$ is even, then $m-1$ is odd and 
by induction hypothesis we have
$$(\Lambda_a-1)^m = (\Lambda_a-1)^{m-1} (\Lambda_a-1)
\equiv  (\Lambda_a-1)^2 \equiv 1 \pmod{2},$$
which shows the assertion for $m$. If $m$ is odd, then
$m-1$ is even and we have
$$(\Lambda_a-1)^m = (\Lambda_a-1)^{m-1} (\Lambda_a-1)
\equiv \Lambda_a-1 \pmod{2},$$
which proves the assertion for $m$.

(2) For $m=1$, the assertion is obvious.
Suppose $m \geq 2$ and that the assertion holds
for $m-1$. Then, we have
$$(\Lambda_a-1)^m = (\Lambda_a-1)^{m-1} (\Lambda_a-1)
\equiv (\Lambda_a-1)^2 =  a\Lambda_a - 2 \Lambda_a + 1
\equiv \Lambda_a - 1 \pmod{2},$$
since $a$ is odd. This proves the assertion for $m$.

(3) We have, by Lemma~\ref{lemma:basic},
$$(\Lambda_a-1)(\Lambda_{ab} -1) = a\Lambda_{ab}
- \Lambda_a - \Lambda_{ab} + 1 =(a-1)\Lambda_{ab}
- \Lambda_a + 1 \equiv \Lambda_a-1 \pmod{2},$$
since $a$ is odd, which shows the assertion.

(4) When $m=1$, the assertion is obvious.
Suppose $m \geq 2$ and that the assertion holds
for $m-1$. Then, we have, by Lemma~\ref{lemma:basic},
\begin{eqnarray*}
\prod_{j=1}^m (\Lambda_{a_j} - 1)
& \equiv & \left(\sum_{j=1}^{m-1} \Lambda_{a_j} - 1
\right) (\Lambda_{a_m} - 1)  \pmod{2} \\
& = & \sum_{j=1}^{m-1}\Lambda_{a_j}\Lambda_{a_m}
- \Lambda_{a_m} - \sum_{j=1}^{m-1} \Lambda_{a_j} + 1 \\
& = & \sum_{j=1}^{m-1}(a_j, a_m)\Lambda_{[a_j, a_m]}
- \sum_{j=1}^m \Lambda_{a_j} + 1 \\
& \equiv & \sum_{j=1}^m \Lambda_{a_j} - 1 \pmod{2},
\end{eqnarray*}
since $a_j$ are all even, which proves the assertion.

(5) We have, by Lemma~\ref{lemma:basic},
$$\Lambda_a (\Lambda_b-1) = (a, b)\Lambda_{[a, b]}
- \Lambda_a \equiv \Lambda_{[a, b]} - \Lambda_a
\pmod{2},$$
since $(a, b)$ is odd, which proves the assertion.

This completes the proof.
\end{proof}

Then, we have the following.

\begin{lemma}\label{lemma:FMe}
We have
$$\prod_{a \in E_f}
(\Lambda_a - 1) \equiv \prod_{a \in \overline{E}_f}
(\Lambda_a - 1)  \pmod{2}.$$
\end{lemma}

\begin{rmk}\label{rmk:zero}
When $\overline{E}_f = \emptyset$, 
$$\prod_{a \in \overline{E}_f}
(\Lambda_a - 1)$$
is understood to be equal to $1$ in the
group ring $\Z \C^*$. 
This is because we have $\overline{E}_f = \emptyset$
if and only if all the exponents are even and each such number
appears an even number of times, and in such a case, by
Lemma~\ref{lemma:FM} (1), the result is $1$ modulo $2$.
\end{rmk}

\begin{proof}[Proof of Lemma~\textup{\ref{lemma:FMe}}]
By Lemma~\ref{lemma:FM} (1) for $m$ even,
even if we perform the procedure as in
Definition~\ref{dfn:Ebar} (\ref{(1)}),
the modulo $2$ class of the product of $\Lambda_a-1$ over all
elements $a$ of the relevant set
does not change.
Then, by Lemma~\ref{lemma:FM} (1) for $m$ odd and (2),
the same holds with the procedure of Definition~\ref{dfn:Ebar} (\ref{(2)}).
Finally, by Lemma~\ref{lemma:FM} (3), after
the procedure as in Definition~\ref{dfn:Ebar} (\ref{(3)}),
we see that
$$
\prod_{a \in E_f} (\Lambda_a - 1) 
\equiv \prod_{a \in \overline{E}_f}(\Lambda_a-1) \pmod{2}.
$$
This completes the proof.
\end{proof}

\begin{proof}[Proof of Theorem~\textup{\ref{thm:FM}}]
We will prove the equivalence based on Brieskorn's
formula (\ref{eq:BR}).

Suppose that $\overline{E}_f = \overline{E}_g$ holds. Then,
by Lemma~\ref{lemma:FMe} together with Lemma~\ref{rem:square},
we see that $\Delta_f(t)$ and $\Delta_g(t)$ satisfy
the Fox--Milnor type relation.

Conversely, suppose that $\Delta_f(t)$
and $\Delta_g(t)$ satisfy the Fox--Milnor type
relation. Then, again by Lemmas~\ref{lemma:FMe}
and \ref{rem:square}, we have
$$\prod_{a \in \overline{E}_f}(\Lambda_a-1) \equiv
\prod_{b \in \overline{E}_g}(\Lambda_b-1) \pmod{2}.
$$
Let $\overline{E}_f^0$ (resp.\ $\overline{E}_f^1$)
be the subset of $\overline{E}_f$ consisting of even (resp.\ odd)
integers. We also define $\overline{E}_g^0$ and 
$\overline{E}_g^1$
similarly. Then, we have
\begin{eqnarray*}
& & \left(\prod_{a \in \overline{E}^0_f}(\Lambda_a-1) \right)
\left(\prod_{a \in \overline{E}^1_f}(\Lambda_a-1) \right) \\
& \equiv & 
\left(\prod_{b \in \overline{E}^0_g}(\Lambda_b-1) \right)
\left(\prod_{b \in \overline{E}^1_g}(\Lambda_b-1) \right)
\pmod{2}.
\end{eqnarray*}
By Lemma~\ref{lemma:FM} (4), we have
\begin{eqnarray}
& & \left(\sum_{a \in \overline{E}^0_f}\Lambda_a-1 \right)
\left(\prod_{a \in \overline{E}^1_f}(\Lambda_a-1) \right) 
\label{eq:FM2}\\
& \equiv & 
\left(\sum_{b \in \overline{E}^0_g}\Lambda_b-1 \right)
\left(\prod_{b \in \overline{E}^1_g}(\Lambda_b-1) \right)
\pmod{2}. \nonumber
\end{eqnarray}
By comparing the terms of the forms $\Lambda_d$ with $d$ odd
on both sides of the above congruence,
with the help of the basic formula (\ref{eq:basic}) together with
the linear independence of $\{\Lambda_a\}_a$ modulo $2$
as described in \cite[Lemma~3.3]{BS2}, we have
$$\prod_{a \in \overline{E}^1_f}(\Lambda_a-1) 
\equiv \prod_{b \in \overline{E}^1_g}(\Lambda_b-1) 
\pmod{2}.
$$
As no integer in $\overline{E}^1_f$ (or
$\overline{E}_g^1$) is a multiple of another one,
by the same argument as in the proof of
\cite[Theorem~2.7]{BS2}, we see that
\begin{equation}
\overline{E}^1_f = \overline{E}^1_g.
\label{eq:E1}
\end{equation}

By (\ref{eq:FM2}), we have
\begin{eqnarray}
& & \left(\sum_{a \in \overline{E}^0_f}\Lambda_a \right)
\left(\prod_{a \in \overline{E}^1_f}(\Lambda_a-1) \right) 
\label{eq:EE}
\\
& \equiv & 
\left(\sum_{b \in \overline{E}^0_g}\Lambda_b \right)
\left(\prod_{b \in \overline{E}^1_g}(\Lambda_b-1) \right)
\pmod{2}. \nonumber
\end{eqnarray}
Then, by considering the terms of the forms
$\Lambda_d$ with $d$ minimal on both sides, we see that
$$\min \overline{E}_f^0 = \min \overline{E}_g^0,$$
which we set as $m_0$.
Consequently, by subtracting 
$$\Lambda_{m_0}\left(\prod_{a \in \overline{E}^1_f}(\Lambda_a-1) \right) 
\equiv 
\Lambda_{m_0}\left(\prod_{b \in \overline{E}^1_g}(\Lambda_b-1) \right) 
\pmod{2}
$$
from both sides of (\ref{eq:EE}), we get
\begin{eqnarray*}
& & \left(\sum_{a \in \overline{E}^0_f \setminus \{m_0\}}
\Lambda_a \right)
\left(\prod_{a \in \overline{E}^1_f}(\Lambda_a-1) \right) 
\\
& \equiv & 
\left(\sum_{b \in \overline{E}^0_g \setminus \{m_0\}}
\Lambda_b \right)
\left(\prod_{b \in \overline{E}^1_g}(\Lambda_b-1) \right)
\pmod{2}. 
\end{eqnarray*}
Repeating this procedure, we finally get
$\overline{E}^0_f = \overline{E}^0_g$.
This together with (\ref{eq:E1}) implies 
$\overline{E}_f = \overline{E}_g$.
This completes the proof
of Theorem~\textup{\ref{thm:FM}}.
\end{proof}

\begin{rmk}\label{rmk:FM}
By \cite[Proposition~2.6]{BS2}, if the Seifert forms
of $K_f$ and $K_g$ are Witt equivalent over the real
numbers (i.e., if they have the same equivariant signatures),
then their Alexander polynomials satisfy the Fox--Milnor
type relation. So, by Theorem~\ref{thm:FM}, we have
$\overline{E}_f = \overline{E}_g$.
\end{rmk}

\begin{corollary}\label{cor1}
Suppose that the exponents of a Brieskorn
polynomial $f$ are all distinct
and that no exponent is a multiple of another odd exponent.
Let $g$ be an arbitrary Brieskorn polynomial
with the same number of variables as $f$.
Then $K_f$ and $K_g$ are cobordant if and only if
they have the same set of exponents.
In particular, if the exponents of $f$ are all even and all distinct,
the same conclusion holds.
\end{corollary}

\begin{proof}
Under the assumption for $f$, we see that
$E_f = \overline{E}_f$ holds by the definition
of the essential exponent set. Suppose that $K_f$ and $K_g$
are cobordant. Then, their Alexander polynomials
satisfy the Fox--Milnor type relation, and by Theorem~\ref{thm:FM},
we have $\overline{E}_f = \overline{E}_g$.
As $E_f = \overline{E}_f$ has $n+1$ distinct elements,
so does $\overline{E}_g$. As this is a subset of $E_g$, we must
have $\overline{E}_g = E_g$. Hence we have $E_f = E_g$.
This completes the proof.
\end{proof}

We also have the following.

\begin{corollary}\label{cor2}
Let $f$ be a Brieskorn polynomial. Then,
the Alexander polynomial $\Delta_f(t)$ of the
algebraic knot $K_f$ associated with $f$ is 
never a square.
\end{corollary}

\begin{proof}
First, 
by the same argument as in the proof of Lemma~\ref{rem:square},
we see that $\Delta_f(t)$ is a square if and only if
$\mathrm{divisor}\,\Delta_f(t) \equiv 0 \pmod{2}$.

If all the exponents are even and each of them appears an even
number of times, then we
have $\overline{E}_f = \emptyset$, and
$\Delta_f(t)$ is not a square
by Lemma~\ref{lemma:FMe} and Remark~\ref{rmk:zero}.

On the other hand, if $E_f$ contains an odd integer, then it
persists in $\overline{E}_f$. Furthermore, if there are no
odd exponents and an even
integer appears exactly an odd number of times, then
it persists in $\overline{E}_f$.
Hence, in these cases, 
$\mathrm{divisor}\,\Delta_f(t) \not\equiv 0 \pmod{2}$
by Lemma~\ref{lemma:FMe}. Therefore,
$\Delta_f(t)$ is not a square.

This completes the proof.
\end{proof}

\begin{corollary}\label{prop:nullcob}
The algebraic knot $K_f$ associated with a
Brieskorn polynomial $f$ is never null-cobordant.
\end{corollary}

\begin{proof}
Suppose that $K_f$ is null-cobordant.
Then its Alexander polynomial $\Delta_f(t)$
and that of the trivial knot satisfy the Fox--Milnor type relation.
As the Alexander polynomial of the trivial
knot is equal to $1$, by 
Lemma~\ref{rem:square}, we see that
$\Delta_f(t)$ must be a square. This contradicts
Corollary~\ref{cor2}. This completes the proof.
\end{proof}

In fact, we have a stronger result as follows.

\begin{proposition}\label{prop:order}
Let $K_f$ be the algebraic knot associated with a
Brieskorn polynomial $f$.
If it is spherical, then it always has infinite order
in the knot cobordism group.
\end{proposition}

\begin{proof}
Suppose $K_f$ is of finite order. Then, its
equivariant signatures all vanish, since they are
integer-valued additive invariants of the cobordism
classes. Therefore, by
Remark~\ref{rmk:FM}, its Alexander polynomial
and that of the trivial knot satisfies the Fox--Milnor type relation.
Hence, by Lemma~\ref{rem:square}, we see that
$\Delta_f(t)$ must be a square. 
This contradicts
Corollary~\ref{cor2}. This completes the proof.
\end{proof}

\begin{rmk}
In fact, Michel \cite{Michel} proves the
following results for a general holomorphic
function germ $f : (\C^{n+1}, \mathbf{0}) \to (\C, 0)$
possibly with an isolated critical point at $\mathbf{0}$.

(1) If the algebraic knot $K_f$ is spherical, then
it is null-cobordant if and only if $f$
does not have a singularity at $\mathbf{0}$.

(2) If $f$ has a singularity at the origin
and $K_f$ is spherical,
then $K_f$ has infinite order 
in the knot cobordism group $C_{2n-1}$.

Our results show that although the Fox--Milnor type
relation seems to be weak, it leads to the above 
important results
in the case of Brieskorn polynomials.
\end{rmk}

\begin{proposition}\label{thm:finite}
Let $K_f$ and $K_g$ be the algebraic knots
associated with Brieskorn polynomials $f$ and
$g$, respectively. We assume
that they are spherical. 
\begin{enumerate}
\item If $K_f \sharp (-K_g^!)$
is of finite order in the knot cobordism group,
then the order must be equal to $1$ or $2$.
\item If $K_f$ and $K_g$ have the same equivariant
signatures, then $K_f \sharp K_f$ is cobordant to
$K_g \sharp K_g$.
\end{enumerate}
\end{proposition}

\begin{proof}
(1) It is known that $K_f \sharp (-K_g^!)$
is of finite order if and only if its equivariant signatures
all vanish. (For example, see \cite{Collins, L2}.)
Therefore, by our assumption,
the equivariant signatures of $K_f$ and $K_g$
coincide, and by Remark~\ref{rmk:FM},
the Alexander polynomials of $K_f$ and $K_g$ satisfy
the Fox--Milnor type relation. Hence,
$\Delta_f(t)\Delta_g(t)$ is a square.
Then, by \cite[Theorem~3.4.8]{Collins}, 
$K_f \sharp (-K_g^!)$ cannot have order $4$
in the knot cobordism group, since its
Alexander polynomial coincides with
$\Delta_f(t)\Delta_g(t)$, which is a square.
Hence, the order must be equal to $1$ or $2$.

(2) Since $K_f \sharp (-K_g^!)$ has order $1$ or $2$,
we see that 
$$2 (K_f \sharp (-K_g^!)) =
(K_f \sharp K_f) \sharp (-(K_g \sharp K_g)^!)$$ 
is null-cobordant,
and the result follows.
\end{proof}

We also have the following.

\begin{proposition}\label{prop:share}
Let $f$ and $g$ be Brieskorn polynomials.
If the algebraic knots $K_f$ and $K_g$ 
are cobordant, then their Alexander
polynomials $\Delta_f(t)$ and $\Delta_g(t)$
share at least one irreducible cyclotomic polynomial factor.
\end{proposition}

Before proving the above proposition, let us recall the 
following result due to Brieskorn \cite{Brieskorn}.
For a Brieskorn polynomial $f$ with the exponent set $E_f$,
we construct a finite graph $\Gamma_f$ as follows: the vertices
correspond to the elements of $E_f$, and for $a, b \in E_f$,
we connect them by an edge if their greatest common divisor
satisfies
$(a, b) > 1$. A connected component of  $\Gamma_f$
is called an \emph{odd $2$--component} if its vertex set consists 
of an odd  number of even integers such that each pair of 
vertices are connected by an edge
and their greatest common divisor is always equal to $2$.
Then we have the following.

\begin{theorem}[Brieskorn \cite{Brieskorn}]\label{thm:Br}
Let $f$ be a Brieskorn polynomial of $n+1$ variables.
For $n \neq 2$, the algebraic knot
$K_f$ is spherical if and only if $\Gamma_f$ satisfies one of the following.
\begin{enumerate}
\item The graph $\Gamma_f$ contains at least two isolated vertices.
\item The graph $\Gamma_f$ contains one isolated vertex and an odd
$2$--component.
\end{enumerate}
\end{theorem}

\begin{rmk}\label{rem:Mumford}
In the above theorem,
for $n=2$, we have $H_*(K_f; \Z) \cong
H_*(S^3; \Z)$ if and only if the above conditions (1) and
(2) hold, which is implicit in the proof due to Brieskorn
\cite{Brieskorn}. On the other hand, by a result of
Mumford \cite{Mum}, $K_f$ is homeomorphic to $S^3$
if and only of $f^{-1}(0)$ does not have a singularity
at the origin.
\end{rmk}

\begin{proof}[Proof of Proposition~\textup{\ref{prop:share}}]
By Theorem~\ref{thm:Br},
we see that by adding appropriate powers of extra two variables
to $f$, we get a Brieskorn polynomial $\tilde{f}$ of $n+3$
variables such that $K_{\tilde{f}}$ is spherical.

Suppose that the equivariant signatures for $K_f$
all vanish. Then, its Seifert form is Witt equivalent to $0$
over the real numbers (see \cite[\S4]{Sa00}).
Since the Seifert form for $K_{\tilde{f}}$
is the tensor product of that
for $K_f$ and a certain matrix (see \cite{Saka}), 
we see that it is also Witt equivalent to $0$
over the real numbers. Hence, its
equivariant signatures all vanish. Since
$K_{\tilde{f}}$ is spherical, this implies that it
has finite order in the knot cobordism group,
which contradicts
Proposition~\ref{prop:order}.
Hence, an equivariant signature of $K_f$ with respect to
a root $\lambda$ of $\Delta_f(t)$ does not
vanish. As an equivariant signature is a cobordism
invariant, the equivariant signature of $K_g$
with respect to $\lambda$ does not vanish, either.
This implies that $\lambda$ is a root of $\Delta_g(t)$.
As the Alexander polynomials $\Delta_f(t)$
and $\Delta_g(t)$ are products of cyclotomic
polynomials, the result follows.
\end{proof}

\section{Linear independence in the knot cobordism 
group}\label{section3.5}

Litherland \cite{Litherland} has shown that
the algebraic knots associated with the Brieskorn
polynomials $z_1^p + z_2^q$ with $2 \leq p < q$
and $(p, q) = 1$ are linearly independent
in the knot cobordism group of dimension $1$.
Note that these algebraic knots in $S^3$ are the so-called
\emph{torus knots}.

In order to prove a similar result in higher dimensions,
let us prepare the following. For a fixed integer $n \geq 1$,
let $\mathcal{B}$ be a set of exponent sets
of $n+1$ elements such that for each exponent set
belonging to $\mathcal{B}$, the exponents are
greater than or equal to $2$ and are pairwise
relatively prime, and that
no two of the exponent
sets of $\mathcal{B}$ have equal product.
In other words, for $\{p_i\}_{i=1}^{n+1} \neq
\{q_i\}_{i=1}^{n+1} \in \mathcal{B}$, we have
$$p_1p_2 \cdots p_{n+1} \neq
q_1q_2 \cdots q_{n+1}.$$
We call such a set $\mathcal{B}$ a \emph{good family of
exponent sets}.
For example, the set $\mathcal{P}$ of all exponent
sets such that the exponents are distinct
prime numbers is a good family of exponent sets.

\begin{theorem}\label{thm:LI}
Let $\mathcal{B}$ be a good family of exponent sets
of $n+1$ elements, and consider
the family of Brieskorn polynomials whose exponent sets correspond
bijectively to the elements of $\mathcal{B}$.
Then, for $n \neq 2$, the associated algebraic knots are spherical
and are linearly independent in the knot cobordism group
of dimension $2n-1$.
\end{theorem}

Note that the corresponding algebraic knots
are easily seen to be spherical by Theorem~\ref{thm:Br}.

For the proof of Theorem~\ref{thm:LI}, 
let us prepare some materials.
Let $K$ be a spherical $(2n-1)$--knot and $L$
its Seifert matrix. For a complex number $\zeta$ of modulus $1$,
let us consider the signature of the Hermitian matrix
$$(1 - \zeta) L + (1 - \bar{\zeta}) L^T.$$
This is independent of the choice of Seifert matrix $L$.
This gives rise to a function $S^1 \to \Z$, where $S^1$
is the unit circle in $\C$, and it is known to be
continuous (and therefore constant) everywhere except at
$(-1)^{n+1}$ times the unit roots of the Alexander
polynomial $\Delta_K(t)$ (for example, 
see \cite[Chapter~9]{Collins}).
This function is not a cobordism invariant in general:
however, the jumps at $(-1)^{n+1}$ times
the unit roots of the Alexander polynomial
are cobordism invariants.
This is called the \emph{signature jump function}.
It is known that it can be written in terms of the equivariant
signatures
(see \cite{L1, Mat} or \cite[Theorem~3.4.7]{Collins}).

Now let $\{p_i\} = \{p_1, p_2, \ldots, p_{n+1}\}$ be 
an exponent set in $\mathcal{B}$ and
set 
$$P = p_1p_2 \cdots p_{n+1}.$$
Note that the integers $p_1, p_2, \ldots, p_{n+1}$ are relatively
prime to each other.
For a positive integer $r$, set
\begin{align*}
L_+\left(\frac{r}{P}\right) 
& = \left\{(k_1, k_2, \ldots, k_{n+1}) \in \Z^{n+1}\,\left|\,
\sum_{i=1}^{n+1} \frac{k_i}{p_i} 
\equiv \frac{r}{P} \pmod{2}, \right.\right.\\
& \quad \left.\, 0 < k_i < p_i, \,
i = 1, 2, \ldots, n+1\bigg\},\right. \\
L_-\left(\frac{r}{P}\right) & 
= \left\{(k_1, k_2, \ldots, k_{n+1}) \in \Z^{n+1}\,\left|\,
\sum_{i=1}^{n+1} \frac{k_i}{p_i} 
\equiv \frac{r}{P} + 1 \pmod{2}, \right.\right.\\
& \quad \left.\, 0< k_i < p_i, \,
i = 1, 2, \ldots, n+1\bigg\}, \right.
\end{align*}
where for rational numbers $s$ and $t$, and a positive
integer $u$, we write
$s \equiv t \pmod{u}$ if the difference $s-r$ is an integer
multiple of $u$.

Then, we have the following.

\begin{lemma}
The set $L_+(r/P) \cup L_-(r/P)$
contains at most one element, and we have that
$L_+(r/P) = L_-(r/P) = \emptyset$ if and only if $r$ is a
multiple of some $p_i$.
\end{lemma}

\begin{proof}
Since $p_1, p_2, \ldots, p_{n+1}$ are pairwise
relatively prime, the following three are equivalent to
each other:
\begin{eqnarray*}
(1) & & \sum_{i=1}^{n+1} \frac{k_i}{p_i} 
\equiv \frac{r}{P} \pmod{1}, \\
(2) & &  \sum_{i=1}^{n+1} k_i p_1 p_2
\cdots \widecheck{p_i} \cdots p_{n+1}
\equiv r \pmod{P}, \\
(3) & & k_i p_1 p_2 \cdots \widecheck{p_i} \cdots p_{n+1}
\equiv r \pmod{p_i}, \quad 1 \leq \forall i \leq n+1,
\end{eqnarray*}
where ``$\widecheck{p_i}$'' means that it is deleted from the
product. Since for each $i$, $p_i$ and 
$p_1 p_2 \cdots \widecheck{p_i} \cdots p_{n+1}$ are relatively
prime, if $r \not\equiv 0 \pmod{p_i}$,
then there exists a unique $k_i \in \{1, 2, \ldots, p_i-1\}$
satisfying the above property (3) for that $i$.
On the other hand, if $r \equiv 0 \pmod{p_i}$ for some $i$,
there does not exist $k_i \in \{1, 2, \ldots, p_i-1\}$ with
the above property (3) for that $i$. These
observations lead to the desired conclusion.
This completes the proof.
\end{proof}

Furthermore, the following follows from \cite[\S9.3]{Collins}.

\begin{lemma}\label{lemma:jump}
The signature jump at $\exp{(2\pi \sqrt{-1} r/P)}$ is equal to $1$
if $L_+(r/P) \neq \emptyset$, 
is equal to $-1$ if $L_-(r/P) \neq \emptyset$,
and is equal to $0$ if $r$ is a multiple of some $p_i$.
\end{lemma}

\begin{proof}[Proof of Theorem~\textup{\ref{thm:LI}}]
Let us show that the signature jump functions $j_{\{p_i\}}$
for the algebraic knots corresponding to the exponent sets $\{p_i\}$
are linearly independent over $\Z$ for $\{p_i\} \in
\mathcal{B}$.
Suppose there is a nontrivial dependence relation
among $j_{\{p_i\}}$ over $\Z$.
Let $M$ be the maximum of $P = p_1p_2 \cdots p_{n+1}$
appearing in a nontrivial dependence relation.
Note that by the definition of a good family of exponent sets,
such maximum is attained only by a unique element
$\{q_i\}$ in $\mathcal{B}$.
Since $j_{\{p_i\}}(1/M) = 
0$ for $\{p_i\}$ with $p_1 p_2 \cdots p_{n+1} < M$,
we see that  the dependence relation implies that
$j_{\{q_i\}}(1/M) = 0$. This contradicts Lemma~\ref{lemma:jump}.
Therefore, the jump functions corresponding to the
elements of $\mathcal{B}$ are linearly independent over $\Z$.
Since the jump functions are additive cobordism invariants,
the result follows.
\end{proof}

\begin{rmk}
The above proof is based on the idea used in \cite{Litherland}
for $n=1$. In Theorem~\ref{thm:LI}, we imposed the
condition that no two of the exponent
sets of $\mathcal{B}$ have equal product. We do not
know if this condition is redundant or not.
\end{rmk}

\section{Decomposition of Seifert form}\label{section4}

As has been noted in \S\ref{section2},
for $n > 1$, we have that the knot cobordism group
$C_{2n-1}$ is isomorphic to $\Z^\infty \oplus \Z_2^\infty 
\oplus \Z_4^\infty$.
Clearly such a decomposition is not unique:
however, according to \cite{L1, L2}, there is a family of
explicitly defined cobordism invariants
that induce such a decomposition canonically
(for details, see \cite[\S3.2 and \S3.4]{Collins}, for example).
For instance, for a representative of a given element of $C_{2n-1}$,
for each irreducible symmetric
polynomial $\delta(t) \in \Q[t]$ that
is a factor of the Alexander polynomial, we have an 
associated element in
the Witt group $G_\Q^\delta$
of isometric structures
over $\Q$ corresponding to a power of $\delta$.
It is also known that
$G_{\Q} \cong \oplus_\delta G_{\Q}^\delta$,
where the sum is over all irreducible symmetric polynomials
$\delta$.

In this section, we show that for some elements
in $C_{2n-1}$ with $n=3$ associated with Brieskorn
polynomials 
$$g(z_1, z_2, z_3, z_4) = z_1^3 + z_2^4 + z_3^4 + z_4^p$$
with $p \geq 5$ prime to $2$ and $3$,
the associated factor in $G_\Q^\delta$
can be trivial for some nontrivial factor $\delta(t)$
of the Alexander polynomial, although the associated 
algebraic knots are of
infinite order in the knot cobordism group.

For $a \geq 2$, let $M_a$ be the $(a-1) \times (a-1)$
unimodular matrix
$$M_a = \begin{pmatrix} 1 & -1 & 0 & \cdots & 0 \\
                                       0 & 1 & -1 & \ddots & \vdots \\
                                     0 & 0 & 1 & \ddots & 0 \\
                                       \vdots   &   \vdots   &   \ddots   &   \ddots  & -1 \\
                                        0 & 0 & \cdots & 0 & 1
                                        \end{pmatrix}.
                                        $$ 
Note that the Seifert form of the algebraic knot
associated with the Brieskorn polynomial
$$f = z_1^{a_1} + z_2^{a_2} + \cdots + z_{n+1}^{a_{n+1}}$$
is given by the tensor product $L = M_{a_1} \otimes 
M_{a_2} \otimes
\cdots \otimes M_{a_{n+1}}$, see \cite{Saka}, in which the
matrix appears to be different from the above one: 
this is due to the difference
in the definition of the Seifert form.
 
Note that we have
\begin{eqnarray}
& & S = L + (-1)^n \overline{L^T}, \quad
H = (-1)^{n+1}\overline{L^{-1}}L^T, 
\label{eq:matrix}
\\
& & T = (-1)^{n+1}L \overline{(L^{-1})^T}, \quad
S = L(I - \overline{H}), \nonumber
\end{eqnarray}
where $L$ is the Seifert matrix, $S$ is the sesquilinearized
intersection form of the Milnor fiber, $H$ is the homological
monodromy matrix, $T$ is the cohomological monodromy matrix,
and $I$ is the identity matrix (for example, see \cite{Sa00}).
Here, we have given the above formula in such a way that
they work also over the complex numbers.

Let us consider an explicit example:
$$f(z_1, z_2, z_3) = z_1^3 + z_2^4 + z_3^4.$$
The divisor of its characteristic polynomial of the monodromy
(or Alexander polynomial) $\Delta_f(t)$ is given by
$$\mathrm{divisor}\, \Delta_f = (\Lambda_3 - 1)
(\Lambda_4 - 1)(\Lambda_4 - 1) = 2 \Lambda_{12}
+ \Lambda_3 - 2 \Lambda_4 - 1,$$
and hence we have
\begin{eqnarray*}
\Delta_f(t) & = & \frac{(t^{12}-1)^2 (t^3-1)}{(t^4-1)^2(t-1)} \\
& = & \frac{\phi_{12}^2 \, \phi_6^2 \, \phi_4^2 \, \phi_3^2 
\, \phi_2^2 \,
\phi_1^2 \, \phi_3 \, \phi_1}{\phi_4^2 \, 
\phi_2^2 \, \phi_1^2 \,\phi_1} \\
& = & \phi_{12}^2 \, \phi_6^2 \, \phi_3^3,
\end{eqnarray*}
where for a positive integer $m$, $\phi_m(t)$ denotes
the $m$--th cyclotomic polynomial.
Note that the degrees of $\phi_{12}, \phi_6, \phi_3$
are equal to $4, 2, 2$, respectively.
According to Steenbrink's formula \cite{Steenbrink2},
the equivariant signatures
corresponding to the roots of
$\phi_{12}, \phi_6, \phi_3$ are equal to
$8, 0, 6$, respectively. (This means, for example, that
the sum of the equivariant signatures corresponding to the 
roots of $\phi_{12}$ is equal to $8$.)
Therefore, the $\phi_{12}$-- and $\phi_3$--primary
components of the Seifert form of $f$ are not zero
in the Witt groups $G_\Q^{\phi_{12}}$
and $G_\Q^{\phi_3}$, respectively,
in the sense of \cite[\S10]{L2} or \cite{Collins}.

Let us analyze the $\phi_6$--primary component of the
Seifert form of $f$ in the Witt group $G_\Q^{\phi_6}$.
The Seifert form for $f$ is given by the
unimodular $(18 \times 18)$--matrix
$$L = M_3 \otimes M_4 \otimes M_4.$$
The form $M_3$ is irreducible over $\Q$, since its
Alexander polynomial $\phi_3$ is irreducible.
On the other hand, the Alexander polynomial of $M_4$
is equal to $\phi_4 \, \phi_2$, which is not irreducible.
Let us decompose $M_4$ into the irreducible factors over $\Q$.

By some computations, we see the following:
\begin{eqnarray*}
M_4 & = & \begin{pmatrix} 1 & -1 & 0 \\
0 & 1 & -1 \\
0 & 0 & 1 \end{pmatrix},\\
M_4^{-1} & = &  \begin{pmatrix} 1 & 1 & 1 \\
0 & 1 & 1 \\
0 & 0 & 1 \end{pmatrix},\\
S_4 & = & M_4 + M_4^T = \begin{pmatrix} 2 & -1 & 0 \\
-1 & 2 & -1 \\
0 & -1 & 2 \end{pmatrix},\\
T_4 & = & - M_4(M_4^{-1})^T  = \begin{pmatrix} 0 & 1 & 0 \\
0 & 0 & 1 \\
-1 & -1 & -1 \end{pmatrix},\\
H_4 & = & (T_4)^T = \begin{pmatrix} 0 & 0 & -1 \\
1 & 0 & -1 \\
0 & 1 & -1 \end{pmatrix}.\\
\end{eqnarray*}
Note that
\begin{eqnarray*}
(H_4)^T M_4 H_4 & = &  \begin{pmatrix} 0 & 1 & 0 \\
0 & 0 & 1 \\
-1 & -1 & -1 \end{pmatrix}
\begin{pmatrix} 1 & -1 & 0 \\
0 & 1 & -1 \\
0 & 0 & 1 \end{pmatrix}
\begin{pmatrix} 0 & 0 & -1 \\
1 & 0 & -1 \\
0 & 1 & -1 \end{pmatrix} \\
& = & \begin{pmatrix} 1 & -1 & 0 \\
0 & 1 & -1 \\
0 & 0 & 1 \end{pmatrix} = M_4.
\end{eqnarray*}
The eigenvalues of $H_4$ are $-1, \pm \sqrt{-1}$.
Eigenvectors corresponding to the eigenvalues
$-1, \sqrt{-1}$ and $-\sqrt{-1}$ are given by
$$\begin{pmatrix}1 \\ 0 \\ 1 \end{pmatrix}, \,
\begin{pmatrix}1 \\ 1 - \sqrt{-1} \\ -\sqrt{-1} \end{pmatrix}, \,
\begin{pmatrix}1 \\ 1 + \sqrt{-1} \\ \sqrt{-1} \end{pmatrix},
$$
respectively. Therefore, the $\phi_2$--primary component
is generated by
$$\begin{pmatrix}1 \\ 0 \\ 1 \end{pmatrix},$$
and the $\phi_4$--primary component is generated by
$$\begin{pmatrix}1 \\ 1 \\ 0 \end{pmatrix}, \,
\begin{pmatrix}0 \\ 1 \\ 1 \end{pmatrix}.
$$
(For this, consider the real and the imaginary parts of
the corresponding eigenvectors.)
Set
$$P = \begin{pmatrix} 1 & 1 & 0 \\
0 & 1 & 1 \\
1 & 0 & 1 \end{pmatrix}.
$$
Then, we have
$$P^{-1} = \frac{1}{2}
\begin{pmatrix} 1 & -1 & 1 \\
1 & 1 & -1 \\
-1 & 1 & 1 \end{pmatrix}
$$
and
$$P^{-1} H_4 P = \begin{pmatrix} -1 & 0 & 0 \\
0 & 0 & -1 \\
0 & 1 & 0 \end{pmatrix}.
$$
So, we have verified that $P$ gives the correct
decomposition of the monodromy into the irreducible
components.

\begin{rmk}
We can show that we cannot choose an integral unimodular
matrix as $P$ as follows. If we choose
$$a\begin{pmatrix}1 \\ 1 \\ 0 \end{pmatrix}
+ b\begin{pmatrix}0 \\ 1 \\ 1 \end{pmatrix} \mbox{ and }
a'\begin{pmatrix}1 \\ 1 \\ 0 \end{pmatrix}
+ b'\begin{pmatrix}0 \\ 1 \\ 1 \end{pmatrix}
$$
as bases for the $\phi_4$--primary component for
some integers $a, b, a', b'$, then
we can show that the determinant of the
corresponding $(3 \times 3)$--matrix is an even integer.
\end{rmk}

Then, we have
$$P^T M_4 P = \begin{pmatrix} 2 & 0 & 0 \\
0 & 1 & -1 \\
0 & 1 & 1 \end{pmatrix}.
$$
So, over $\Q$, the bilinear form $M_4$ is
isomorphic to $(1) \oplus R$, where
$$R = \begin{pmatrix}
1 & -1 \\ 1 & 1 \end{pmatrix}.
$$
Then, we have, over $\Q$,
\begin{eqnarray*}
L & = & M_3 \otimes M_4 \otimes M_4 \\
& \cong & M_3 \otimes ((1) \oplus R) \otimes ((1) \oplus
R) \\
& \cong & M_3 \otimes ((1) \oplus R \oplus R \oplus (R \otimes R)) \\
& \cong & M_3 \oplus (M_3 \otimes R) \oplus
(M_3 \otimes R) \oplus (M_3 \otimes R \otimes R).
\end{eqnarray*}
The characteristic polynomials corresponding to the
four irreducible factors are given by
$$\phi_3(t), \phi_{12}(t), \phi_{12}(t) \mbox{ and }
\phi_3(t)^2\phi_6(t)^2,$$
respectively. So, in order to analyze the $\phi_6$--primary component
of $L$, we still need to decompose $M_3 \otimes R \otimes R$,
which is an $(8 \times 8)$--matrix.

Note that by (\ref{eq:matrix}), the monodromy matrix 
$H_3$ corresponding to
$M_3$ is given by
$$H_3 = - M_3^{-1}M_3 =  \begin{pmatrix}
0 & -1 \\ 1 & -1 \end{pmatrix}.
$$
We have
$$R^{-1} = \frac{1}{2}\begin{pmatrix}
1 & 1 \\ -1 & 1 \end{pmatrix}
$$
and the monodromy matrix $H_R$ corresponding to $R$
is given by
$$H_R = - R^{-1}R^T = \begin{pmatrix}
0 & -1 \\ 1 & 0 \end{pmatrix}.
$$
Set $\omega = \exp{(2 \pi \sqrt{-1}/3)}$.
Eigenvectors of $H_3$ corresponding to the eigenvalues 
$\omega$ and $\bar{\omega}$ are given by
$$u_1 = \begin{pmatrix}
1 \\ -\omega \end{pmatrix} \mbox{ and }
u_2 = \begin{pmatrix}
1 \\ -\bar{\omega} \end{pmatrix},
$$
respectively. Eigenvectors of $H_R$ corresponding to
the eigenvalues $\sqrt{-1}$ and $-\sqrt{-1}$ are given by
$$v_1 = \begin{pmatrix}
1 \\ -\sqrt{-1} \end{pmatrix} \mbox{ and }
v_2 = \begin{pmatrix}
1 \\ \sqrt{-1} \end{pmatrix},
$$
respectively. Therefore, the monodromy matrix $H_{3, R, R}$
associated with
$M_3 \otimes R \otimes R$ is diagonalized by the
$(8 \times 8)$--matrix $Q$ consisting of the $8$ column vectors
$$u_i \otimes v_j \otimes v_k,$$
$i, j, k = 1, 2$, in such a way that
$$Q^{-1}  H_{3, R, R} Q = \begin{pmatrix}
\omega & 0 \\
0 & \bar{\omega}
\end{pmatrix} \otimes
\begin{pmatrix}
\sqrt{-1} & 0 \\
0 & -\sqrt{-1}
\end{pmatrix}
\otimes
\begin{pmatrix}
\sqrt{-1} & 0 \\
0 & -\sqrt{-1}
\end{pmatrix}.
$$
Therefore, the $\phi_6$--primary component is generated by
$u_1 \otimes v_1 \otimes v_1$,
$u_1 \otimes v_2 \otimes v_2$,
$u_2 \otimes v_1 \otimes v_1$ and
$u_2 \otimes v_2 \otimes v_2$
over $\C$.
Note that
\begin{eqnarray*}
u_1 \otimes v_1 \otimes v_1
& = & (1, -\sqrt{-1}, -\sqrt{-1}, -1, -\omega, \omega\sqrt{-1},
\omega\sqrt{-1}, \omega)^T,\\
u_1 \otimes v_2 \otimes v_2
& = &  (1, \sqrt{-1}, \sqrt{-1}, -1, -\omega, -\omega\sqrt{-1},
-\omega\sqrt{-1}, \omega)^T,\\
u_2 \otimes v_1 \otimes v_1
& = & (1, -\sqrt{-1}, -\sqrt{-1}, -1, -\bar{\omega}, 
\bar{\omega}\sqrt{-1},
\bar{\omega}\sqrt{-1}, \bar{\omega})^T,\\
u_2 \otimes v_2 \otimes v_2
& = & (1, \sqrt{-1}, \sqrt{-1}, -1, -\bar{\omega}, 
-\bar{\omega}\sqrt{-1},
-\bar{\omega}\sqrt{-1}, \bar{\omega})^T.
\end{eqnarray*}
By considering the real and imaginary parts, we see that 
the $\phi_6$--primary component is generated by
\begin{eqnarray*}
w_1 & = &  (1, 0, 0, -1, 1/2, -\sqrt{3}/2,
-\sqrt{3}/2, -1/2)^T,\\
 w_2 & = & (0, -1, -1, 0, -\sqrt{3}/2, -1/2,
-1/2, \sqrt{3}/2)^T,\\
 w_3 & = & (1, 0, 0, -1, 1/2, \sqrt{3}/2,
\sqrt{3}/2, -1/2)^T,\\
 w_4 & = & (0, 1, 1, 0, -\sqrt{3}/2, 1/2,
1/2, \sqrt{3}/2)^T
\end{eqnarray*}
over $\R$. Then, we have
\begin{eqnarray*}
w_1 + w_3 & = &  (2, 0, 0, -2, 1, 0,
0, -1)^T,\\
(w_1 - w_3)/\sqrt{3} & = & (0, 0, 0, 0, 0, -1,
-1, 0)^T,\\
(w_2 + w_4)/\sqrt{3} & = & (0, 0, 0, 0, -1, 0,
0, 1)^T,\\
-(w_2 - w_4) & = & (0, 2, 2, 0, 0, 1,
-, 0)^T.
\end{eqnarray*}
Note that these four vectors can be written as
\begin{eqnarray*}
r_1 & = & (2, 1)^T \otimes (1, 0, 0, -1)^T, \\
r_2 & = & (0, -1)^T \otimes (0, 1, 1, 0)^T, \\
r_3 & = & (0, -1)^T \otimes (1, 0, 0, -1)^T, \\
r_4 & = & (2, 1)^T \otimes (0, 1, 1, 0)^T,
\end{eqnarray*}
respectively. Then, by calculating
$$r_i^T (M_3 \otimes R \otimes R) r_j,$$
$i, j = 1, 2, 3, 4$, we see that the $\phi_6$--primary
component of the bilinear form $M_3 \otimes R \otimes R$
is isomorphic over $\Q$ to
$$\begin{pmatrix} 0 & -4 & 0 & -12 \\
-4 & 0 & 4 & 0 \\
0 & -4 & 0 & 4 \\
12 & 0 & 4 & 0
\end{pmatrix},
$$
which, in turn, is isomorphic to
$$\begin{pmatrix} 0 & 0 & -1 & -3 \\
0 & 0 & -1 & 1 \\
-1 & 1 & 0 & 0 \\
3 & 1 & 0 & 0
\end{pmatrix}
$$
over $\Q$.

Consequently, we see that the $\phi_6$--primary component
of $M_3 \otimes R \otimes R$, hence that of $L$,
is algebraically null-cobordant, since it admits a metabolizer
(see \cite[\S3.1]{Collins}, for example).

Let $p$ be a positive integer
relatively prime to $2$ and $3$, and consider the
Brieskorn polynomial
$$g(z_1, z_2, z_3, z_4) = z_1^3 + z_2^4 + z_3^4
+ z_4^p.$$
Then, the algebraic knot $K_g$ associated with $g$
is spherical, i.e.\ $K_g$ is homeomorphic to the sphere
$S^7$ (see Theorem~\ref{thm:Br}).
As the Seifert matrix $L_g$ of $K_g$ is given by
the tensor product of $L$ and $M_p$,
we see that a certain direct summand of $L_g$
is null-cobordant over $\Q$ in the algebraic cobordism group.

\medskip

The above explicit example shows that if we consider 
the image of the
cobordism class of a spherical algebraic knot associated with
a Brieskorn polynomial in $\oplus_\delta G_{\Q}^\delta$,
then there might be a direct summand which vanishes in 
$G_{\Q}^\delta$ for some $\delta$.
This means that even if two algebraic knots
are cobordant, the irreducible factors of their Alexander
polynomials might be different, although they share
at least one irreducible factor according to Proposition~\ref{prop:share}.

This suggests a major difficulty in proving the
cobordism invariance of the exponents for Brieskorn polynomials.
Note that in \cite{YS}, topological invariance of
exponents for Brieskorn polynomials was proved using
the topological invariance of the Alexander polynomial.
Such an approach seems not to work for the study of
cobordisms.

\section{Cyclic suspension}\label{section5}

In this section, we explore cyclic suspensions
of simple fibered knots and algebraic knots, and also 
their properties concerning
cobordisms.

Let $K \subset S^{2n+1}$ be a $(2n-1)$--knot.
Then, we can move the standard sphere
$S^{2n+1} \subset S^{2n+3}$ ambient isotopically
to get $S'$ such that $S'$ intersects $S^{2n+1}$
transversely along $K$. For a positive integer $d$,
we consider the $d$--fold cyclic branched covering 
$\widetilde{S}$ of $S^{2n+3}$
branched along $S^{2n+1}$, which is diffeomorphic to
$S^{2n+3}$. Then the pull-back $K_d$ 
of $S'$ by the branched covering
map in $\widetilde{S}$ is called the \emph{$d$--fold
cyclic suspension} of $K$. Furthermore, we call the positive
integer $d$ the \emph{suspension degree}.
Note that $K_d$ itself
is diffeomorphic to the $d$--fold cyclic branched 
covering of $S^{2n+1}$
branched along $K$, and that it is considered to be
a $(2n+1)$--knot in $S^{2n+3}$.
This notion has been introduced by Kauffman \cite{Kauffman} 
and Neumann \cite{N} (see also \cite{KN}).
Note that if $K$ is a simple fibered knot,
then so is $K_d$.

In this section, we consider the following problem.

\begin{prob}\label{prob:sus}
For a common integer $d$,
let $(K_i)_d$ be the $d$--fold
cyclic suspensions of two knots $K_i$, $i = 1, 2$.
Furthermore, for another common integer $e$,
let $(K_i)_{d, e}$ be the $e$--fold cyclic
suspensions of $(K_i)_d$, $i = 1, 2$.
Is it possible to construct examples such that
$K_i$ are not cobordant, that
$(K_i)_d$ are cobordant and that
$(K_i)_{d, e}$ are not cobordant?
\end{prob}

If the answer is affirmative, then it would show that
the cyclic suspensions do not preserve cobordisms in general.

Note that the algebraic knot associated with a
Brieskorn polynomial $z_1^{a_1} + z_2^{a_2} +
\cdots + z_{n+1}^{a_{n+1}}$ is the iterated
cyclic suspension of the $(a_1, a_2)$--torus link in $S^3$
(see \cite{N}).
The above problem is closely related to the study of
cobordisms of such knots.

Let $n \geq 3$ be an integer.
For the moment, we will assume that $n$ is odd.
Consider the matrices
$$A_1 = \begin{pmatrix}
B & C \\
-C^T & \mathbf{0}
\end{pmatrix}$$
and
$$A_2 = \begin{pmatrix} 0 & 1 \\ -1 & 0
\end{pmatrix},
$$
where $B$ is a $2 \times 2$ integer matrix
with $\det(B + B^T) = \pm 1$, $C$
is a $2 \times 2$ integer matrix with $\det C = \pm 1$,
and $\mathbf{0}$ denotes the $2 \times 2$ zero matrix.
So, $A_1$ is a unimodular $(4 \times 4)$--matrix and
$A_2$ is a unimodular $(2 \times 2)$--matrix.
Let $K_1$ and $K_2$ be the simple fibered
$(2n-1)$--knots in $S^{2n+1}$ whose Seifert matrices are given by
$A_1$ and $A_2$, respectively. Such simple fibered knots
exist uniquely up to isotopy by \cite{Durfee, Kato}.

For positive integers $a$ and $b$,
let $(K_i)_a$ be the $a$--fold cyclic suspension of the
knot $K_i$, and $(K_i)_{a, b}$
be the $b$--fold cyclic suspension of $(K_i)_a$, $i = 1, 2$.
Then, their Seifert matrices $(A_i)_a$ and $(A_i)_{a, b}$,
respectively, are given by
$$(A_i)_a = A_i \otimes M_a \mbox{ and }
(A_i)_{a, b} = A_i \otimes M_a \otimes M_b$$
(see \cite{Kauffman, KN, Neumann}).

Let us consider the $2$--fold cyclic suspensions
$(K_i)_2$, $i = 1, 2$. 
As $M_2$ is the $(1 \times 1)$--matrix $(1)$,
we can identify their Seifert matrices with those of $K_i$,
$i = 1, 2$.
As we have
$$(S_1)_2 = A_1 + A_1^T = \begin{pmatrix}
B + B^T & \mathbf{0} \\
\mathbf{0} & \mathbf{0} \end{pmatrix}, \quad
(S_2)_2 = A_2 + A_2^T = \mathbf{0},$$
we see that 
$$H_n((K_1)_2; \Z) \cong H_n((K_2)_2; \Z)
\cong \Z \oplus \Z \cong H_{n-1}((K_1)_2; \Z) 
\cong H_{n-1}((K_2)_2; \Z).$$
(For example, see the argument just after \cite[Remark~5.9]{BS}.)
Furthermore, as $A_1$ and $A_2$ both have metabolizers,
so does $A_1 \oplus (-A_2)$.

However, 
$(K_1)_2$ and $(K_2)_2$ are not cobordant,
since the Seifert forms restricted to
$H_n((K_i)_2; \Z) = \Ker S_i$, $i = 1, 2$, are not isomorphic
(see \cite{BM}).
Note that these knots are not spherical.
 
Let us now consider the $3$--fold cyclic suspensions
$(K_1)_3$ and $(K_2)_3$, respectively.
Then, their Seifert matrices are given by
$$(A_1)_3 = A_1 \otimes M_3
= \begin{pmatrix}
B \otimes M_3 & C \otimes M_3 \\
-C^T \otimes M_3 & \mathbf{0}
\end{pmatrix}
$$
and
$$(A_2)_3 = A_2 \otimes M_3 =
\begin{pmatrix} \mathbf{0} & M_3 \\
-M_3 & \mathbf{0}
\end{pmatrix},
$$
respectively. Then, the intersection matrices of their fibers are given
by
\begin{eqnarray*}
(S_1)_3 & = & (A_1)_3 + (A_1)_3^T
= \begin{pmatrix}
B \otimes M_3 + B^T \otimes M_3^T & C \otimes M_3 
-C \otimes M_3^T \\ C^T \otimes M_3^T - C^T \otimes M_3
& \mathbf{0}
\end{pmatrix} \\
& = & \begin{pmatrix}
B \otimes M_3 + B^T \otimes M_3^T & C \otimes (M_3 
- M_3^T) \\ - C^T \otimes (M_3- M_3^T)
& \mathbf{0}
\end{pmatrix}
\end{eqnarray*}
and
$$(S_2)_3 = (A_2)_3 + (A_2)_3^T
= \begin{pmatrix}
\mathbf{0} & M_3 - M_3^T  \\ 
-(M_3 - M_3^T)  & \mathbf{0}
\end{pmatrix},
$$
respectively.
As we have
$$\det (M_3 - M_3^T) = \det
\begin{pmatrix} 0 & -1 \\ 1 & 0 \end{pmatrix}
= 1,
$$
we see that  both $(S_1)_3$ and $(S_2)_3$ are
unimodular. Therefore, the fibered knots
$(K_1)_3$ and $(K_2)_3$ are spherical
(for example, see the argument just after \cite[Remark~5.9]{BS}).
As their Seifert matrices are obviously algebraically null-cobordant,
the knots are, in fact, null-cobordant, and in particular
they are cobordant.
 
We can also show that $K_1$ and $K_2$ are not
diffeomorphic to each other for an appropriate choice of $C$.
For example, consider
$$C  = \begin{pmatrix} 0 & 1 \\ -1 & 0 \end{pmatrix}.$$
In this case, the intersection matrices are
$$S_1 = A_1 - A_1^T = \begin{pmatrix}
B - B^T & C + C^T \\
-(C + C^T) & \mathbf{0} \end{pmatrix}
= \begin{pmatrix} B - B^T & \mathbf{0} \\
\mathbf{0} & \mathbf{0}
\end{pmatrix}
$$
and
$$S_2 = A_2 - A_2^T = \begin{pmatrix}
0 & 2 \\ 2 & 0 \end{pmatrix}.
$$
Therefore, the rank of $H_{n-1}(K_1; \Z)$
is greater than or equal to $2$, while $H_{n-1}(K_2; \Z)$
is finite of order $4$.
So, $K_1$ and $K_2$ are not diffeomorphic and hence
are not cobordant.

Let us now consider $(K_1)_{2, 3} = (K_1)_{3, 2}$
and $(K_2)_{2, 3} = (K_2)_{3, 2}$.
Their Seifert forms are given by
$$(A_1)_{3, 2} = A_1 \otimes M_3 \otimes M_2 =
\begin{pmatrix}
B \otimes M_3 & C \otimes M_3 \\
-C^T \otimes M_3 & \mathbf{0}
\end{pmatrix}
$$
and
$$(A_2)_{3, 2} = A_2 \otimes M_3 \otimes M_2 
 =
\begin{pmatrix} \mathbf{0} & M_3 \\
-M_3 & \mathbf{0}
\end{pmatrix},
$$
respectively. 
Then, their intersection matrices are
$$(S_1)_{3, 2} = (A_1)_{3, 2} - (A_1)_{3, 2}^T
= \begin{pmatrix}
B \otimes M_3 - B^T \otimes M_3^T & C \otimes
M_3 + C \otimes M_3^T \\
-C^T \otimes M_3 - C^T \otimes M_3^T & \mathbf{0}
\end{pmatrix}
$$
and
$$(S_2)_{3, 2} = (A_2)_{3, 2} - (A_2)_{3, 2}^T
= \begin{pmatrix}
\mathbf{0} & M_3 + M_3^T \\
-M_3 - M_3^T & \mathbf{0}
\end{pmatrix},
$$
respectively.
For $C$ as above, we see that
$$|\det (S_1)_{3, 2}|  = 3^4, \quad
|\det (S_2)_{3, 2} | = 3^2,$$
and hence $(K_1)_{3, 2}$ and $(K_2)_{3, 2}$
are not diffeomorphic to each other and are not cobordant.
 
Summarizing, we have the following.
\begin{enumerate}
\item $K_1$ and $K_2$ are not diffeomorphic to each other and 
are not cobordant.
\item $(K_1)_2$ and $(K_2)_2$ are diffeomorphic to each other, 
but are not cobordant.
\item $(K_1)_3$ and $(K_2)_3$ are spherical and 
null-cobordant, so
they are cobordant to each other.
\item $(K_1)_{3, 2} = (K_1)_{2, 3}$ and $(K_2)_{3, 2}
= (K_2)_{2, 3}$ are not diffeomorphic to each other
and are not cobordant.
\end{enumerate}
 
So, this answers Problem~\ref{prob:sus}
affirmatively.

\begin{rmk}
In general, 
if $K_1$ and $K_2$ are spherical knots
which are cobordant, then $(K_1)_{2, 2}$ and
$(K_2)_{2, 2}$ are also cobordant. See \cite[\S8]{KN}.
\end{rmk}

Now, let us consider examples of algebraic knots.
In \cite{dBM}, Du Bois and Michel constructed two polynomials 
$$f = h_{r, s, p, q}(z_1, z_2, \ldots, z_{n+1}) \mbox{ and }
g = h_{s-8, r+8, p, q}(z_1, z_2, \ldots, z_{n+1})$$
with isolated critical points at the origin such that
$K_f$ and $K_g$ are cobordant, although they are not isotopic.
Let $k$ be a positive integer called an \emph{exponent}
in the sense of \cite{dBM} for both of $f$ and $g$:
i.e., $(t_f^k -1)^2$ and $(t_g^k - 1)^2$ both vanish, where
$t_f$ and $t_g$ are homological monodromies for the
Milnor fibrations for $f$ and $g$,
respectively, and ``$1$'' denotes the identity homomorphism.
Let us consider the algebraic knots $K_{\tilde{f}}$
and $K_{\tilde{g}}$ associated with
$$\tilde{f}(z_1, z_2, \ldots, z_{n+2})
= f(z_1, z_2, \ldots, z_{n+1}) + z_{n+2}^k$$
and
$$\tilde{g}(z_1, z_2, \ldots, z_{n+2})
= g(z_1, z_2, \ldots, z_{n+1}) + z_{n+2}^k,$$
respectively.
Note that they are $k$--fold cyclic suspensions of $K_f$ and $K_g$,
respectively.

\begin{lemma}
The homology groups
$H_n(K_{\tilde{f}}; \Z)$ and
$H_n(K_{\tilde{g}}; \Z)$ have non-isomorphic torsions.
\end{lemma}

\begin{proof}
Recall that $K_{\tilde{f}}$ (resp.\ $K_{\tilde{g}}$)
is the $k$--fold cyclic branched cover of $S^{2n+1}$
branched along $K_f$ (resp.\ $K_g$).
This implies that $K_{\tilde{f}}$ admits an open book
structure with page diffeomorphic to $F_f$ and with
algebraic monodromy $t = t_f^k$.

Let $B \subset K_{\tilde{f}}$ be the branched locus
and let $E$ be the complement of an open tubular
neighborhood of $B$ in $K_{\tilde{f}}$.
Thus, $E$ is the total space of a fiber bundle
over $S^1$ with fiber $F_f$ and with algebraic monodromy
$t = t_f^k$. Then, we have the following Wang exact sequence
of homology \cite{Wang} (see also \cite[Lemma~8.4]{Milnor}):
$$H_n(F_f; \Z) \spmapright{t-1} H_n(F_f; \Z)
\to H_n(E; \Z) \to H_{n-1}(F_f; \Z).$$
Since $F_f$ is $(n-1)$--connected \cite{Milnor},
we have $H_{n-1}(F_f; \Z) = 0$ so that  we have
$$H_n(E; \Z) \cong H_n(F_f; \Z)/\Image(t - 1).$$
Then, by the Meyer--Vietoris exact sequence
for the pair $(E, N(B))$, where $N(B)$ is the
closed tubular neighborhood of $B$ in $K_{\tilde{f}}$,
we have that
$$H_n(\partial N(B); \Z) \to H_n(N(B); \Z) \oplus
H_n(E; \Z) \to H_n(K_{\tilde{f}}; \Z) \to H_{n-1}(\partial N(B); \Z)$$
is exact. As $N(B) \cong K_f \times D^2$ and $K_f$ is homeomorphic
to $S^{2n-1}$ with $n \geq 3$,
we see that $H_n(\partial N(B); \Z)$, $H_n(N(B); \Z)$ and
$H_{n-1}(\partial N(B); \Z)$ all vanish.
Therefore, we have
$H_n(K_{\tilde{f}}; \Z) \cong H_n(E; \Z)$,
and hence they are isomorphic to the quotient
$H_n(F_f; \Z)/(t_f^k - 1)H_n(F_f; \Z)$.

On the other hand, $\Ker(t_f^k-1)$ is a pure submodule
of the free abelian group $H_n(F_f; \Z)$ of finite rank,
where a submodule of a free abelian group is said to be
\emph{pure} if it is a direct summand.
Therefore, there exists a free abelian subgroup $H_f$
of $H_n(F_f; \Z)$ such that
$H_n(F_f; \Z) = H_f \oplus \Ker(t_f^k-1)$.
As $\Image(t_f^k -1)$ is contained in $\Ker(t_f^k-1)$,
we see that $H_n(K_{\tilde{f}}; \Z) \cong
H_n(F_f; \Z)/(t_f^k - 1)H_n(F_f; \Z)$
is isomorphic to $H_f \oplus 
\left(\Ker(t_f^k-1)/\Image(t_f^k-1)\right)$.
Note that a similar isomorphism holds for $H_n(K_{\tilde{g}}; \Z)$
as well.

Since the twist groups, which are the torsion
subgroups of  $\Ker(t_f^k-1)/\Image(t_f^k-1)$
and $\Ker(t_g^k-1)/\Image(t_g^k-1)$, are not isomorphic
to each other
according to \cite{dBM}, we see that the torsion subgroups
of $H_n(K_{\tilde{f}}; \Z)$ and 
$H_n(K_{\tilde{g}}; \Z)$ are not isomorphic.
\end{proof}

The above lemma implies that although $K_f$ and
$K_g$ are cobordant, their cyclic suspensions
$K_{\tilde{f}}$
and $K_{\tilde{g}}$ are not, since they are not diffeomorphic.

If we take further iterated cyclic suspensions
appropriately, say $K_{\hat{f}}$ and $K_{\hat{g}}$, where
$$\hat{f}(z_1, z_2, \ldots, z_{n+3}, z_{n+4})
= \tilde{f}(z_1, z_2, \ldots, z_{n+2}) + z_{n+3}^v
+ z_{n+4}^w$$ and
$$\hat{g}(z_1, z_2, \ldots, z_{n+3}, z_{n+4})
= \tilde{g}(z_1, z_2, \ldots, z_{n+2}) + z_{n+3}^v
+ z_{n+4}^w$$ for some appropriate prime numbers $v$ and $w$,
then $K_{\hat{f}}$ and $K_{\hat{g}}$ are
spherical and hence are cobordant.

Summarizing, we have the following.
\begin{enumerate}
\item The algebraic knots $K_f$ and $K_g$ are
cobordant, but are not isotopic.
\item Their $k$--fold cyclic suspensions $K_{\tilde{f}}$
and $K_{\tilde{g}}$ are not diffeomorphic to each other 
and are not cobordant.
\item The iterated cyclic suspensions
$K_{\hat{f}}$ and $K_{\hat{g}}$ of $K_{\tilde{f}}$
and $K_{\tilde{g}}$, respectively, are cobordant.
\end{enumerate}

This is yet another example that shows that
cyclic suspensions (with a fixed suspension degree)
do not behave well with respect
to cobordisms. This time, the example shows this phenomenon
for algebraic knots.

\begin{rmk}\label{rem:spherical}
(1) If $K_0$ and $K_1$ are cobordant knots, then
if their cyclic suspensions $\tilde{K}_0$ and
$\tilde{K}_1$, respectively, of the same degree are spherical
of dimension greater than or equal to $3$,
then they are cobordant. This is because the Seifert
matrices of $\tilde{K}_i$ are tensor products of those
of $K_i$, which are (algebraically) cobordant, 
and the same matrix, and hence they are algebraically
cobordant. For spherical higher dimensional knots,
this implies that they are cobordant (see \cite{L1}).

(2) Similarly, if $K$ is a spherical knot which has
finite order in the knot cobordism group, then
if its cyclic suspension $\tilde{K}$ is spherical,
then $\tilde{K}$ also has finite order in the
knot cobordism group. This is because, since
the Seifert form of $K$ is Witt equivalent to $0$
over the real numbers, so is that of $\tilde{K}$.

(3) Suppose that $K$ is a spherical knot and that
its $d$--fold cyclic suspension $\tilde{K}$ is also spherical.
Let us suppose that $\tilde{K}$ is null-cobordant.
Then, we do not know if $K$ is also null-cobordant or not.

Similarly, suppose that $K_0$ and $K_1$
are spherical knots and that their $d$--fold
cyclic suspensions $\tilde{K}_0$ and $\tilde{K}_1$,
respectively, are also spherical. Let us suppose that
$\tilde{K}_0$ and $\tilde{K}_1$ are cobordant.
Then, we do not know if $K_0$ and $K_1$ are also
cobordant or not, except for the case $d=2$.
\end{rmk}

Since the algebraic knots associated with Brieskorn polynomials are 
iterated cyclic suspensions of torus knots, the observations
in this section 
may show that by adding extra variables we may encounter 
a pair of algebraic knots associated with Brieskorn polynomials 
which are cobordant but which have distinct exponents. 

\section*{Acknowledgment}\label{ack}
The authors would like to thank Professor Mutsuo Oka
for valuable discussions.
The second author would like to thank 
the Institut de Recherche Math\'{e}matique Avanc\'{e}e 
(IRMA, UMR 7501) of CNRS and the University of Strasbourg
for the hospitality
during the preparation of the manuscript.
This work has been supported in part by JSPS KAKENHI 
Grant Numbers 
JP22K18267, JP23H05437.


\begin{thebibliography}{99999}
%
\bibitem{BM}V.~Blanl\oe il and F.~Michel,  
{\em A theory of cobordism for non-spherical links}, 
Comment.\ Math.\ Helv.\ \textbf{72} (1997), 30--51.  
%
\bibitem{BS}V.~Blanl\oe il and O.~Saeki, 
{\em Cobordism of fibered knots and related topics},
in ^^ ^^ Singularities in geometry and topology 2004", 
pp.~1--47, Adv.\ Stud.\ Pure Math.\ \textbf{46}, 
Math.\ Soc.\ Japan, Tokyo, 2007.
%
\bibitem{BS2}V.~Blanl\oe il and O.~Saeki, 
{\em Cobordism of algebraic knots defined by Brieskorn polynomials},
Tokyo J.\ Math.\ \textbf{34} (2011), 429--443.
%
\bibitem{Brieskorn}
E.~Brieskorn, {\em Beispiele zur Differentialtopologie von 
Singularit\"{a}ten},
Invent.\ Math.\ \textbf{2} (1966), 1--14.
%
\bibitem{Brieskorn70}
E.~Brieskorn, {\em Die Monodromie der isolierten Singularit\"{a}ten
von Hyperfl\"{a}chen},
Manuscripta Math.\ \textbf{2} (1970), 103--161.

\bibitem{Collins}J.~Collins, {\em On the concordance 
orders of knots},
Thesis, Doctor of Philosophy,
University of Edinburgh, 2012,
arXiv:1206.0669 [math.GT].
%
\bibitem{dBM}P.~Du Bois and F.~Michel,
{\em Cobordism of algebraic knots via Seifert forms},
Invent.\ Math.\ \textbf{111} (1993), 151--169.
%
\bibitem{Durfee}A.~Durfee, {\em Fibered knots and algebraic 
singularities}, Topology \textbf{13} (1974), 47--59.
%
\bibitem{F-M0}R.~H.~Fox and J.~W.~Milnor, 
{\em Singularities of $2$--spheres in 
$4$--space and equivalence of knots}, Bull.\ Amer.\ 
Math.\ Soc.\ \textbf{63} (1957), 406.
%
\bibitem{F-M}R.~H.~Fox and J.~W.~Milnor, 
{\em Singularities of $2$--spheres in 
$4$--space and cobordism of knots}, 
Osaka J.\ Math.\ \textbf{3} (1966), 257--267.
%
\bibitem{Kato}M.~Kato, {\em A classification of simple spinnable 
structures on a $1$--connected Alexander manifold}, 
J.\ Math.\ Soc.\ Japan \textbf{26} (1974), 454--463.
%
\bibitem{Kauffman}L.~Kauffman, {\em Branched coverings, 
open books, and knot periodicity}, Topology 
\textbf{13} (1974), 143--160.
%
\bibitem{KN}L.H.~Kauffman and W.D.~Neumann, 
{\em Products of knots, branched fibrations and sums of singularities},
Topology \textbf{16} (1977), 369--393.
%
\bibitem{L1}J.~Levine, 
{\em Knot cobordism groups in 
codimension two}, 
Comment.\ Math.\ Helv.\ \textbf{44} (1969),  
229--244.
%
\bibitem{L2}J.~Levine,
{\em Invariants of knot cobordism},
Invent.\ Math.\ \textbf{8} (1969), 98--110; 
addendum, ibid.\ \textbf{8} (1969), 355.
%
\bibitem{Litherland}R.A.~Litherland, 
{\em Signatures of iterated torus knots},
in ``Topology of low-dimensional manifolds (Proc.\ Second 
Sussex Conf.,
Chelwood Gate, 1977)'', pp.~71--84,
Lecture Notes in Math., 722, Springer, Berlin, 1979.
%
\bibitem{Mat}T.~Matumoto, {\em On the signature invariants of 
a non-singular complex sesquilinear form},
J.\ Math.\ Soc.\ Japan \textbf{29} (1977), 67--71.
%
\bibitem{Michel}F.~Michel, {\em Forme de
Seifert et singularit\'{e}s isol\'ees}, 
in: C.~Weber (ed.), N\oe uds, Tresses et Singularit\'es, 
Plans-sur-Bex 1982, Monogr.\ Enseign.\ Math., vol.~31,
pp.~175--190, Gen\`eve, Enseign.\ Math., 1983.
%
\bibitem{Milnor}J.~Milnor, 
{\sl Singular points of complex hypersurfaces}, 
Ann.\ of Math.\ Stud., Vol.~61, 
Princeton Univ.\ Press, Princeton, N.J.; Univ.\ of Tokyo Press, 
Tokyo, 1968. 
%
\bibitem{MO}J.~Milnor and P.~Orlik, 
{\em Isolated singularities defined by weighted 
homogeneous polynomials},
Topology \textbf{9} (1970), 385--393.
%
\bibitem{Mum}D.~Mumford, {\em The topology of normal
singularities of an algebraic surface and a criterion for simplicity},
Inst.\ Hautes \'{E}tudes Sci.\ Publ.\ Math.\ No.~9 (1961), 
5--22.
%
\bibitem{N}W.D.~Neumann, {\em 
Cyclic suspension of knots and periodicity of signature for singularities},
Bull.\ Amer.\ Math.\ Soc.\ \textbf{80} (1974), 977--981.
%
\bibitem{Neumann}W.D.~Neumann, 
{\em Invariants of plane curve singularities},
N\oe uds, tresses et singularit\'e
(Plans-sur-Bex, 1982), pp.~223--232, 
Monogr.\ Enseign.\ Math., Vol.~31, Enseignement Math., Geneva, 1983. 
%
\bibitem{Sa89}O.~Saeki, 
{\em Topological types of complex isolated hypersurface singularities},
Kodai Math.\ J.\ \textbf{12} (1989), 23--29.
%
\bibitem{Sa00}O.~Saeki, {\em Real Seifert form determines 
the spectrum for semiquasihomogeneous hypersurface 
singularities in $\C^3$}, J.\ Math.\ Soc.\ Japan 
\textbf{52} (2000), 409--431.
%
\bibitem{Saka}K.~Sakamoto, {\em The Seifert matrices of 
Milnor fiberings defined by holomorphic functions},
J.\ Math.\ Soc.\ Japan \textbf{26} (1974), 714--721.
%
\bibitem{SSS}R.~Schrauwen, J.~Steenbrink, and J.~Stevens,
{\em Spectral pairs and the topology of curve singularities},
Complex geometry and Lie theory (Sundance, UT, 1989), pp.~305--328, 
Proc.\ Sympos.\ Pure Math., Vol.~53, Amer.\ Math.\ Soc., 
Providence, RI, 1991. 
%
\bibitem{Steenbrink2}J.H.M.~Steenbrink,
{\em Intersection form
for quasihomogeneous singularities}, Compositio Math.\ 
\textbf{34} (1977), 211--223.

\bibitem{Wang}
H.-C.~Wang, {\em The homology groups of the fibre bundles over 
a sphere},
Duke Math.\ J.\ \textbf{16} (1949), 33--38.

\bibitem{YS}E.~Yoshinaga and M.~Suzuki, {\em On the 
topological types of singularities of Brieskorn--Pham type}, 
Sci.\ Rep.\ Yokohama Nat.\ Univ.\ Sect.\ I \textbf{25} 
(1978), 37--43.

\end{thebibliography}
\end{document}